\documentclass[12pt]{amsart}
\usepackage{amssymb,amsfonts,latexsym,amscd}

\usepackage[all]{xy}
\usepackage{verbatim}
\usepackage{xurl}
\usepackage{hyperref}
\usepackage{mathrsfs}
\usepackage{mathtools}
\usepackage{stmaryrd}
\usepackage{tikz-cd}

\usepackage{xcolor}

\numberwithin{equation}{section}

\newtheorem{theorem}{Theorem}[section]
\newtheorem{lemma}[theorem]{Lemma}
\newtheorem{proposition}[theorem]{Proposition}
\newtheorem{corollary}[theorem]{Corollary}
\theoremstyle{definition}
\newtheorem{definition}[theorem]{Definition}

\newtheorem{example}[theorem]{Example}

\newtheorem{remark}[theorem]{Remark}

\newcommand{\id}{{\rm id}}

\newcommand{\End}{{\rm End}}
\newcommand{\Fun}{{\rm Fun}}

\newcommand{\FPdim}{\text{\rm FPdim}}
\newcommand{\Ext}{\text{\rm Ext}}

\newcommand{\Hom}{{\rm Hom}}
\newcommand{\Ad}{{\rm Ad}}

\newcommand{\Ind}{{\rm Ind}}
\newcommand{\gr}{\text{gr}}

\newcommand{\Rep}{{\rm Rep}}

\newcommand{\Spec}{\text{Spec}}

\newcommand{\tr}{\text{tr}\,}

\newcommand{\eps}{\varepsilon}

\renewcommand{\O}{\mathscr{O}}

\newcommand{\ot}{\otimes}

\newcommand{\ben}{\begin{enumerate}}
\newcommand{\een}{\end{enumerate}}

\newcommand{\C}{{\mathscr C}}
\newcommand{\D}{{\mathscr D}}
\newcommand{\M}{{\mathscr M}}

\newcommand{\E}{\mathcal{E}}

\newcommand{\Z}{{\mathscr{Z}}}

\usepackage{microtype}
\usepackage{booktabs}
\newcommand{\Coh}{\mathrm{Coh}}
\newcommand{\Gm}{\mathbb G_m}
\newcommand{\cR}{\mathcal R}

\newcommand{\Jac}{\mathrm{Jac}}
\newcommand{\soc}{\mathrm{soc}}
\newcommand{\LL}{\mathrm{LL}}

\newcommand{\A}{\mathscr A}
\newcommand{\B}{\mathscr B}
\newcommand{\N}{\mathscr N}
\newcommand{\one}{\mathbf 1}
\hypersetup{
  pdftitle={On finite group scheme-theoretical categories, II},
  pdfauthor={Shlomo Gelaki and Guillermo Sanmarco},
  hidelinks,
  pdfkeywords={finite group schemes, tensor categories, twisted equivariant sheaves}
}

\date{September 29, 2026}

\begin{document}

\title[On finite group scheme-theoretical categories, II]{On finite group scheme-theoretical categories, II}

\author{Shlomo Gelaki}
\address{Department of Mathematics,
Iowa State University, Ames, IA 50011, USA.}
\email{gelaki@iastate.edu}

\author{Guillermo Sanmarco}
\address{Department of Mathematics, University of Washington, Seattle, WA 98195, USA.}
\email{sanmarco@uw.edu}

\makeatletter
\@namedef{subjclassname@2020}{%
	\textup{2020} Mathematics Subject Classification}
\makeatother

\subjclass[2020]{18M20, 16T05, 17B37}
\keywords{finite group schemes; (bi)equivariant sheaves; finite tensor categories;
exact module categories; unimodularity; twisted Drinfeld centers}

\begin{abstract}
Let $\C:=\C(G,\omega,H,\psi)$ be a finite group scheme-theoretical category over an algebraically closed field of characteristic $p\geq0$. For every indecomposable exact module category over $\C$, we classify the simple objects by double cosets and projective representations of their stabilizers. We construct the projective covers using finite-dimensional algebras of twisted equivariant sheaves on the full infinitesimal neighborhoods of the orbits, and determine their underlying sheaves and dimensions. We prove that $\C$ is unimodular if and only if $\delta_G|_H=\delta_H^2$, where $\delta_G$ and $\delta_H$ are the modular characters. In particular, unimodularity is independent of the admissible twisting. We also describe the orbit-indexed Serre summands and relative-center equivariantization of the center of $\Coh(G,\omega)$. In characteristic two, a calculation for $G=\alpha_2$ shows that a nontrivial associator can change both the Loewy length of a projective cover and the dimension of the self-extension space of the simple object.
\end{abstract}

\maketitle
\tableofcontents

\section{Introduction}\label{sec:Introduction}

The main purpose of this paper is to study finite group scheme-theoretical categories
\[
\C:=\C(G,\omega,H,\psi)
\]
with an arbitrary associator, over an algebraically closed field $k$ of characteristic $p\geq0$. Here $G$ is a finite group scheme, $\omega\in Z^3(G,\Gm)$, $H\subseteq G$ is a closed subgroup scheme, and $\psi\in C^2(H,\Gm)$ satisfies $d\psi=\omega|_H$ \cite{G1}. Our results extend the untwisted setting of \cite{GS}, but also address features that are not determined by the classification of simple objects.

The category $\Coh(G,\omega)$ has convolution tensor product and associator $\omega$. Its indecomposable exact module categories are the categories $\M(H,\psi)$ of right $(H,\psi)$-equivariant sheaves \cite[Theorem~5.3]{G1}, with the equivalence relation recalled in Definition~\ref{defequivpairs}; see also \cite[Theorem~1.1]{N} for the finite-group fusion case. The category $\C(G,\omega,H,\psi)$ is the dual with respect to $\M(H,\psi)$, using opposite composition of module endofunctors. Categorical Morita theory \cite[Theorem~7.12.16]{EGNO}, together with Lemma~\ref{newsimple}, identifies its indecomposable exact module categories with
\[
\M((H,\psi),(K,\eta))
=\Coh^{((H,\psi),(K,\eta))}(G,\omega),\qquad d\eta=\omega|_K.
\]
We study the simple objects and projective covers of these categories.

The simple objects are governed by the $(H,K)$-double cosets. If $Z=HgK$ and $L^g=H\cap gKg^{-1}$, the associator and the two cochains determine a cocycle $\xi_g$ on $L^g$; see \eqref{defnrmwg}. Transitive twisted descent gives an equivalence
\[
\Ind_Z:\Rep(L^g,\xi_g)_k\xrightarrow{\simeq}\M_Z,
\qquad
\M=\bigoplus_{Z\in(G/(H\times K))(k)}\overline{\M_Z}.
\]
Here $\M_Z$ consists of objects annihilated by the ideal of the orbit subscheme $Z$, and $\overline{\M_Z}$ is its Serre closure in $\M$. Theorem~\ref{simmodrep} shows that every simple object belongs to exactly one $\M_Z$. The Serre summands need not be indecomposable blocks. In the finite-group fusion case, the orbit description goes back to \cite{O}.

Projective covers require the full infinitesimal neighborhood of an orbit. Let $X_Z\subseteq G$ be the union of the full local components meeting $Z$. We construct a finite-dimensional algebra $A_Z$, with abelian equivalence $A_Z\text{-}{\rm mod}\simeq\overline{\M_Z}$, and a nilpotent ideal $J_Z$ with $(A_Z/J_Z)\text{-}{\rm mod}\simeq\M_Z$. A primitive idempotent for a projective stabilizer representation lifts through $J_Z$ and produces the required projective cover. For a simple $V\in\Rep(L^g,\xi_g)_k$, Theorem~\ref{prsimmodrep} gives
\[
P_\M(\Ind_Z(V))\cong
\O(X_Z)\otimes_k P_{(L^g,\xi_g)}(V)
\]
as underlying $\O(G)$-modules. Formula~\eqref{projective-coaction} recovers the full equivariant structure from $A_Z$. Corollary~\ref{orbit-filtration} describes the orbit-ideal filtration, and Corollary~\ref{semisimplicity} gives a semisimplicity criterion for $\C$.

A second main result is a necessary and sufficient condition for unimodularity. With the integral convention \eqref{eq:modular-convention}, Theorem~\ref{thm:criterion} states that
$\C(G,\omega,H,\psi)$ is unimodular if and only if $\delta_G|_H=\delta_H^2$.
Thus unimodularity depends on the inclusion $H\subseteq G$, but not on the admissible choices of $\omega$ and $\psi$. The proof identifies the equivariant line $\operatorname{Ann}_{\O(X_H)}(I_H)$, where $I_H:=\ker(\O(X_H)\to\O(H))$, and uses it to compute the socle of the projective cover of the unit. For height-one group schemes, Corollary~\ref{cor:lie} gives the equivalent condition
\[
\operatorname{tr}(\operatorname{ad}_{\mathfrak g}x)
=2\operatorname{tr}(\operatorname{ad}_{\mathfrak h}x)
\qquad(x\in\mathfrak h).
\]
Examples~\ref{borel-counterexample} and \ref{modular-cancellation} show that the two modular contributions need not vanish separately. This is a concrete group-scheme criterion alongside the general theory of unimodular module categories in \cite{Y}; related Frobenius conditions for tensor functors and their centralizers are studied in \cite{JY}.

The need to retain the associator in the projective structure can already be seen for $\alpha_2$ in characteristic two. Let $\O(\alpha_2)=k[x]/(x^2)$, with $x$ primitive, and set
$\omega_\lambda:=1+\lambda x\otimes x\otimes x$ $(\lambda\in k)$.
Theorem~\ref{alpha2center} gives an equivalence of abelian categories
\[
\Z(\Coh(\alpha_2,\omega_\lambda))
\simeq k[X,T]/(X^2,T^2-\lambda X)\text{-}{\rm mod}.
\]
For $\lambda=0$ this algebra is $k[X,T]/(X^2,T^2)$, while for $\lambda\ne0$ it is isomorphic to $k[T]/(T^4)$. The Loewy length of the unique projective cover changes from three to four, and the dimension of $\Ext^1(\one,\one)$ changes from two to one (Corollary~\ref{alpha2extensions}). Thus these centers are not equivalent even as abelian categories, although their orbit and simple-object data agree.

The paper is organized as follows. In \S\ref{sec:Preliminaries} we fix the equivariance and comodule conventions, and construct the finite-dimensional algebra used throughout. Section~\ref{S:Biequivariant sheaves on group schemes} proves transitive twisted descent. Section~\ref{sec:doublecosetsandbiequivariant} computes the double-coset multiplier, its gauge covariance and the change of representative. Section~\ref{sec:Module categories over GCTC} treats simple objects, projective covers and fiber functors. In \S\ref{S:structure-GSC} we apply these results to $\C$, including unimodularity and the normal-subgroup obstructions. Finally, \S\ref{sec:twisted-double} treats the center, the $\alpha_2$ example and the relative-center equivariantization, without assuming semisimplicity.

Some assertions in the untwisted setting also require correction. The projective construction replaces the splitting argument, and Examples~\ref{normal-untwisted-counterexample}, \ref{borel-counterexample}, \ref{restriction-not-tensor} and \ref{reduced-inclusion-fails} address the normal-subgroup, unimodularity and center-functor assertions, respectively. Precise comparisons with \cite{GS} are given at these examples.

\medskip
\noindent
{\bf Acknowledgments.} The work of S.G. was supported by Simons Foundation
Award 963288.

\section{Preliminaries}\label{sec:Preliminaries}

\subsection{Conventions}\label{sec:Conventions}
We work over an algebraically closed field $k$ of characteristic $p\ge0$. 

All schemes $X$ considered in this paper are assumed to be {\em finite} over $k$, unless otherwise stated. We specify a point when an augmentation is needed. In particular, a group scheme is pointed by its identity. A general orbit morphism need not preserve a chosen point. Cochains on group schemes are normalized at the identity; functions defining an action-groupoid cocycle are normalized at identity arrows.
For a finite $k$-scheme $X$, we write $|X|:=\dim_k\O(X)$ for its length, and $|X(k)|$ for the number of its closed points. In particular, $|G|$ denotes the order of a finite group scheme $G$, not the order of its group of closed points $G(k)$. We identify the finite set $G(k)$ with its associated constant group scheme when it occurs in a product or a morphism of schemes.
Also, $C^n(G,\Gm)$ consists of normalized invertible regular functions on $G^n$, and $Z^n(G,\Gm)$ is its subgroup of cocycles.

An identity written on points of a group scheme means an identity on points over every commutative $k$-algebra, or, equivalently, an identity of regular functions. Checking only the closed points is not sufficient for a nonreduced scheme.

We use $\alpha_{X,Y,Z}:(X\ot Y)\ot Z\longrightarrow X\ot(Y\ot Z)$ for the associativity constraint. For $\Coh(G,\omega)$ this map acts by $\omega$. The conventions for coalgebras and comodules in \S\ref{sec:Equivariant sheaves and comodules} are derived from this choice. We use the dual convention with opposite composition of module endofunctors, so that $\C(G,\omega,1,1)\simeq\Coh(G,\omega)$.

We assume familiarity with the theory of finite tensor categories and their module categories, and refer to \cite{EGNO} for any unexplained notion.  

\subsection{Sheaves on finite schemes}\label{sec:coh decomposition}
For a finite scheme $X$ as in \S\ref{sec:Conventions}, let ${\rm Coh}(X)$ be the abelian category of coherent sheaves on $X$, i.e., finite-dimensional modules over the finite-dimensional commutative algebra $\O(X)$. Given a closed point $x\in X(k)$, let ${\rm Coh}(X)_x\subseteq {\rm Coh}(X)$ be the abelian subcategory of sheaves on $X$ supported on $x$, so  
\begin{equation*}
{\rm Coh} (X)= \bigoplus_{x \in X(k)}{\rm Coh}(X)_x
\end{equation*}
as abelian categories, where each ${\rm Coh} (X)_x$ contains a unique (up to isomorphism) simple object $\delta_x$ of ${\rm Coh} (X)$. We denote by $P_{x}:=P(\delta_x)$ its projective cover, which also lies in ${\rm Coh} (X)_x$.

Let $Y$ be a finite scheme as in \S\ref{sec:Conventions}, $\varphi:Y\to X$ a scheme morphism, and $\varphi^{\sharp}:\O(X)\to\O(Y)$ the corresponding algebra homomorphism. Recall that $\varphi$ induces a pair $(\varphi^*,\varphi_*)$ of adjoint functors of abelian categories
\begin{equation}\label{invimage}
\varphi^*:{\rm Coh}(X)\to {\rm Coh}(Y),\,\,\,S\mapsto S\ot_{\O(X)}\O(Y),
\end{equation}
where $\O(X)$ acts on $\O(Y)$ via $\varphi^{\sharp}$, and  
\begin{equation}\label{dirimage}
\varphi_*:{\rm Coh}(Y)\to {\rm Coh}(X),\,\,\,T\mapsto T_{\mid \O(X)},
\end{equation}
where $\O(X)$ acts on $T$ via $\varphi^{\sharp}$.

\subsection{Finite group schemes}\label{S:Finite group schemes}
A finite group scheme $G$ is a finite scheme whose coordinate algebra $\O(G)$ is a (finite-dimensional commutative) Hopf algebra (see, e.g., \cite{J,W}). 
It is called {\em \'etale} if $G=G(k)$, where $G(k)$ is the group of closed points of $G$; these finite group schemes are in correspondence with finite abstract groups (see, e.g., \cite[Section 6.4]{W}). On the other hand, $G$ is {\em connected} if $G(k)=1$. For example, 
if $\mathfrak g$ is a finite-dimensional $p$-Lie algebra over $k$ ($p>0$), its $p$-restricted universal enveloping algebra $u(\mathfrak{g})$ is a finite-dimensional cocommutative Hopf algebra (see, e.g., \cite{J}), and   the dual commutative Hopf algebra $u(\mathfrak{g})^*$ is local, and satisfies $x^p=0$ for all $x$ in the augmentation ideal. Thus, the associated group scheme is connected.
More generally, the Frobenius filtration relates finite connected group schemes to group schemes of height one (see \cite{DG,J}).

By a theorem of Cartier (see \cite[Section 11.4]{W}), if $p=0$ then every finite group scheme is \'etale. On the other hand, if $p>0$ then any finite group scheme $G$ is an extension of a connected group scheme by an \'etale one. More precisely, $G$ fits into a split exact sequence
\begin{equation}\label{ses0}
1\to G^{\circ}\xrightarrow{i} G \mathrel{\mathop{\rightleftarrows}^{\pi}_{q}} G(k)\to 1,
\end{equation}
where $G^{\circ}$ is connected and $G(k)$ is \'etale.
In particular, $G(k)$ acts on ${\rm Coh}(G)_1\cong{\rm Coh}(G^{\circ})$, say via $a\mapsto T_a$, and we have an equivalence
\begin{equation}\label{G(k)crossedproduct}
{\rm Coh}(G)\simeq{\rm Coh}(G^{\circ})\rtimes G(k).
\end{equation}
Namely, ${\rm Coh}(G)={\rm Coh}(G^{\circ})\boxtimes {\rm Coh}(G(k))$ as abelian categories, and for any $X_1,X_2\in {\rm Coh}(G^{\circ})$ and $a_1,a_2\in G(k)$, we have
\[\left(X_1\boxtimes a_1 \right)\ot \left(X_2\boxtimes a_2 \right)= \left(X_1\ot T_{a_1}\left(X_2\right)\right)\boxtimes a_1a_2.\]

\subsection{Quotients and free actions}\label{sec:quotmorphind} 
Let $C$ be a finite group scheme, and $X$ a finite scheme as in \S\ref{sec:Conventions}, equipped with a {\em right} $C$-action 
\begin{equation}\label{mu}
\mu:=\mu_{X\times C}:X\times C\to X.
\end{equation}
Equivalently, the algebra homomorphism 
\begin{equation}\label{musharp}
\mu^{\sharp}:\O(X)\to \O(X)\ot \O(C)
\end{equation}
endows $\O(X)$ with a structure of a right $\O(C)$-comodule algebra (see \cite{A}).
Since $C$ is finite, there exists a quotient finite scheme  
\begin{equation}\label{quotient morphismx}
\pi:X\twoheadrightarrow X/C,
\end{equation}
with coordinate algebra 
\begin{equation}\label{doublecosetfa}
\O(X/C)=\O(X)^{C}:=\{f\in\O(X)\mid \mu^{\sharp}(f)=f\ot 1\}.
\end{equation}

Recall that the action $\mu$ (\ref{mu}) is called {\em free} if the morphism 
\begin{equation*}
(p_1,\mu):X\times C\to X\times X
\end{equation*} 
is a closed immersion, where $p_1:X\times C\to X$ is the obvious projection morphism. Recall \cite[Theorem 1(B), p.~112]{Mu} that in this case, the morphism $(p_1,\mu)$ induces a scheme isomorphism 
\begin{equation*}
X\times C\xrightarrow{\cong} X\times_{X/C} X.
\end{equation*}
Equivalently, $(p_1,\mu)$ induces an algebra isomorphism 
\begin{gather*}
\O(X)\otimes_{\O(X/C)} \O(X)\xrightarrow{\cong} \O(X)\otimes_k \O(C),\\
f\ot_{\O(X/C)}\tilde{f}\mapsto (f\ot 1)\mu^{\sharp}(\tilde{f}).
\end{gather*}

\subsection{(Bi)equivariant sheaves}\label{sec:Equivariant sheaves}
Retain the right action $\mu:X\times F\to X$ of a finite group scheme $F$, and write $x\cdot a$ for this action. A normalized action-groupoid $2$-cocycle is an invertible function
\[
\Phi\in\O(X\times F\times F)^{\times},
\]
such that
\begin{equation}\label{groupoid-cocycle}
\begin{split}
\Phi(x;a,b)\Phi(x;ab,c)
&=\Phi(x\cdot a;b,c)\Phi(x;a,bc),\\
\Phi(x;1,a)&=\Phi(x;a,1)=1.
\end{split}
\end{equation}
No group structure on $X$ is required. For a normalized groupoid $1$-cochain $b\in\O(X\times F)^\times$, with $b(x;1)=1$, set
\[
(db)(x;a,a'):=\frac{b(x;a)b(x\cdot a;a')}{b(x;aa')}.
\]

\begin{definition}\label{defequiv}
An $(F,\Phi)$-equivariant sheaf on $X$ is a sheaf $S$ in $\Coh(X)$ together with an isomorphism
\[
\theta:p_1^*S\xrightarrow{\cong}\mu^*S
\]
on $X\times F$, whose restriction to identity arrows is the identity, and which satisfies
\begin{equation}\label{transport-law}
\theta_b(x\cdot a)\theta_a(x)=\Phi(x;a,b)\theta_{ab}(x).
\end{equation}
Here, \eqref{transport-law} denotes an equality between maps of pulled-back sheaves on $X\times F\times F$. 

A morphism between two equivariant sheaves is a map of sheaves commuting with $\theta$. We denote the resulting category by $\Coh^{(F,\Phi)}(X)$.
\end{definition}

\begin{remark}
If $\Phi=1$, this is the usual category $\Coh^F(X)$ of $F$-equivariant sheaves on $X$. \qed
\end{remark}

\begin{example}\label{ex0}
For $\Phi=1$, the structure sheaf $\O(X)$ with its usual $F$-equivariant structure is an object of $\Coh^F(X)$, where the action of $F$ on $\O(X)$ is induced by $\mu^{\sharp}$. \qed
\end{example}

\begin{example}\label{ex1}
Let $X:=\Spec(k)$, and $\xi\in Z^2(F,\Gm)$. Then \eqref{transport-law} says that
$\theta_b\theta_a=\xi(a,b)\theta_{ab}$.
Thus, $\Coh^{(F,\xi)}(\Spec(k))$ is the category of finite-dimensional left comodules over the coalgebra $\O(F)_\xi$, with comultiplication $\Delta_\xi(f)=\Delta(f)\xi$, which we denote by ${\rm Rep}(F,\xi)_k$.

Equivalently, this is the category of right representations of the affine (non-finite) group scheme central extension
\[
F_\xi:=F\times\Gm,\qquad
(a,t)(b,s)=(ab,\xi(a,b)ts),
\]
on which the central $\Gm$ acts by its tautological character. 

For a simple $V\in\Rep(F,\xi)_k$, let $P_{(F,\xi)}(V)$ denote its projective cover. We write $P_F(V)$ when $\xi=1$. \qed
\end{example}

\begin{remark}\label{remarkleftequiv}
A left action can be written as a right action using the inverse. If $B$ acts on $X$ from the left, and $C$ from the right, and the actions commute, we use the right $B\times C$-action
\[
x\cdot(b,c):=b^{-1}xc.
\]
In the twisted setting the two equivariance maps also satisfy a mixed relation. For sheaves on a group scheme this relation is derived from the ambient associator in \S\ref{sec:Equivariant sheaves and comodules}. \qed
\end{remark}

\begin{definition}\label{defequivbi}
Let $G$ be a finite group scheme, $\omega\in Z^3(G,\Gm)$, and let $H,K\subseteq G$ and $\psi,\eta$ satisfy $d\psi=\omega_{\mid H}$ and $d\eta=\omega_{\mid K}$. Set $F:=H\times K$, and use the action $x\cdot(h,k):=h^{-1}xk$. The category
\[
\Coh^{((H,\psi),(K,\eta))}(G,\omega)
\]
is $\Coh^{(F,\Phi)}(G)$, where $\Phi$ is given by \eqref{intro-tau}. A biequivariant sheaf on an invariant subscheme $X\subseteq G$ means the restriction of this same action-groupoid datum to $X$.
\end{definition}

The verification of \eqref{groupoid-cocycle} for \eqref{intro-tau} is given in Proposition \ref{taulemma}. The one-sided notations are
\[
{}^{(H,\psi)}\Coh(G,\omega)
:=\Coh^{((H,\psi),(1,1))}(G,\omega),
\]
\[
\Coh^{(K,\eta)}(G,\omega)
:=\Coh^{((1,1),(K,\eta))}(G,\omega).
\]

\subsection{(Bi)equivariant morphisms}\label{sec:Equivariant morphisms}
Let $\varphi:Y\to X$ be $F$-equivariant, and set
\[
\Phi_Y:=(\varphi\times\id\times\id)^*\Phi_X.
\]
The source and target squares obtained from $\varphi$ and the action maps are cartesian. Thus, the usual base-change identifications for finite morphisms transport equivariance under pullback and pushforward.

\begin{proposition}\label{equivmorphs}
The following hold:
\begin{enumerate}
\item
The adjunction $(\varphi^*,\varphi_*)$ of \eqref{invimage}--\eqref{dirimage} lifts to an adjunction
\[
\Coh^{(F,\Phi_X)}(X)
\mathrel{\mathop{\rightleftarrows}^{\varphi^*}_{\varphi_*}}
\Coh^{(F,\Phi_Y)}(Y).
\]
\item 
The functor $\varphi_*$ is exact. The functor $\varphi^*$ is right exact, and is exact when $\varphi$ is flat. 
\item
If $\varphi$ is a closed immersion, $\varphi_*$ is fully faithful, and its image consists of the objects annihilated by the defining ideal of $Y$.
\end{enumerate}
\end{proposition}

\begin{proof}
Pull back the isomorphism $\theta$ and its identity \eqref{transport-law}. For pushforward use base change in the source and target squares to obtain an isomorphism between the two pullbacks of $\varphi_*S$. Its composition identity is the pushforward of \eqref{transport-law}. The unit and counit of the ordinary adjunction commute with these isomorphisms by naturality. Exactness and the assertion for a closed immersion follow from the corresponding statements for modules over finite-dimensional algebras.
\end{proof}

In particular, Proposition \ref{equivmorphs} applies to one-sided and biequivariant sheaves when the complete twisting datum is pulled back. If $I\subseteq\O(X)$ is $F$-invariant, then $IS\subseteq S$ is an equivariant subobject since both pullbacks of $IS$ are preserved by $\theta$. This observation will be used in the support decomposition.

\begin{example}\label{fromreptocoh}
If $\Phi$ is pulled back from a $2$-cocycle $\xi$ on $F$, the map $X\to\Spec(k)$ gives a functor
\[
{\rm Rep}(F,\xi)_k\longrightarrow\Coh^{(F,\Phi)}(X),
\qquad V\longmapsto\O(X)\ot_k V.
\]
The hypothesis on the whole cocycle, not only on its stabilizer restriction, is needed. \qed
\end{example}

We next realize $\Coh^{(F,\Phi)}(X)$ as the module category of a finite algebra. We use the terminology of corings as in \cite[Chapter~3]{BW}.

Set $R:=\O(X)$ and $\cR:=\O(X\times F)$, with $R$-bimodule structure
\[
r\cdot c\cdot s=p_1^{\sharp}(r)c\mu^{\sharp}(s).
\]
Identify $\cR\ot_R\cR$ with $\O(X\times F\times F)$ by
\[
(c\ot d)(x,a,b)=c(x,a)d(x\cdot a,b),
\]
and define
\begin{equation}\label{coring-comultiplication}
\epsilon(c)(x):=c(x,1),\qquad
\Delta_\Phi(c)(x,a,b):=\Phi(x;a,b)c(x,ab).
\end{equation}
Writing $\Delta_\Phi(c)=\sum c_1\ot_R c_2$, set
\begin{equation}\label{descent-algebra-definition}
A(X,F,\Phi):=\Hom_{R^{\rm op}}(\cR,R),\qquad
(f\star g)(c):=\sum f\bigl(g(c_1)c_2\bigr).
\end{equation}

\begin{proposition}\label{descent-algebra}
The following hold.
\begin{enumerate}
\item The maps in \eqref{coring-comultiplication} make $\cR$ an $R$-coring.
\item Formula~\eqref{descent-algebra-definition} defines an algebra, and there is an equivalence
\[
\Coh^{(F,\Phi)}(X)\simeq A(X,F,\Phi)\text{-}{\rm mod}
\]
preserving underlying $R$-modules.
\item The algebra $A(X,F,\Phi)$ contains $R$ by $r\mapsto r\epsilon$, is free of rank $|F|$ as a left $R$-module, and has dimension $|X||F|$ over $k$.
\end{enumerate}
\end{proposition}

\begin{proof}
Part (1) is the cocycle identity \eqref{groupoid-cocycle} and normalization. A counital left coaction $S\to\cR\ot_R S$ is the same as a transport satisfying \eqref{transport-law}; applying that identity to $a,a^{-1}$ makes the transport invertible.

The change of variables $(x,a)\mapsto(x\cdot a,a^{-1})$ identifies the right $R$-module $\cR$ with a free module of rank $|F|$. For a left coaction $s\mapsto\sum s_{-1}\ot_R s_0$, define
\[
f\cdot s:=\sum f(s_{-1})s_0,
\qquad f\in A(X,F,\Phi),\quad s\in S.
\]
Coassociativity gives $f\cdot(g\cdot s)=(f\star g)\cdot s$, and counitality gives $\epsilon\cdot s=s$. 

Conversely, a right $R$-basis $c_i$ of $\cR$, with dual basis $f_i$, gives the coaction
\begin{equation}\label{dual-basis-coaction}
s\longmapsto\sum_i c_i\ot_R f_i s
\end{equation}
on any left $A(X,F,\Phi)$-module. The multiplication and unit identities give the coaction identities. These constructions are inverse and independent of the basis; the same coassociativity calculation proves associativity of $\star$.

Finally, $r\epsilon$ acts by the original $R$-action, and $r\epsilon=0$ implies $r=0$ since $\epsilon(1)=1$. The equality $(r\epsilon)\star f=rf$ identifies the left $R$-module of the dual algebra with the free dual module of rank $|F|$. This proves (2) and (3).
\end{proof}

\begin{corollary}\label{algebra-base-change}
Set $A:=A(X,F,\Phi)$. If $I\subseteq \O(X)$ is an $F$-invariant ideal, and $X'\subseteq X$ is its closed subscheme, then $IA=AI$, and
\[
A(X,F,\Phi)/IA\cong A(X',F,\Phi_{\mid X'}),
\]
as algebras. In particular, $\Coh^{(F,\Phi)}(X)$ is a finite abelian category, and $IS\subseteq S$ is an equivariant subobject for every $S\in\Coh^{(F,\Phi)}(X)$.
\end{corollary}

\begin{proof}
Invariance and the inverse action show that the ideals generated in $\cR$ by the source and target images of $I$ agree. Then by \eqref{descent-algebra-definition}, $IA$ is a two-sided ideal of $A$. Since $\cR$ is free as a right $R$-module, its dual commutes with reduction modulo $I$. Both the comultiplication and the dual multiplication reduce to the indicated formulas on $X'$. The last statements follow either from this algebra description, or directly from the transport maps.
\end{proof}

\subsection{\texorpdfstring{The tensor category $\Coh(G,\omega)$}{The tensor category Coh(G, omega)}}\label{sec:Coh(G,w)}
Let $G$ be a finite group scheme and let $\omega\in\O(G)^{\ot3}$ be invertible, satisfying
\begin{equation}\label{drinfassoc}
(\id\ot\id\ot\Delta)(\omega)(\Delta\ot\id\ot\id)(\omega)
=(1\ot\omega)(\id\ot\Delta\ot\id)(\omega)(\omega\ot1),
\end{equation}
and
\[(\eps\ot\id\ot\id)(\omega)=(\id\ot\eps\ot\id)(\omega)
=(\id\ot\id\ot\eps)(\omega)=1.
\]
The category $\Coh(G)$ with convolution tensor product and associator
\[
\alpha_{U,V,W}:(U\ot V)\ot W\longrightarrow U\ot(V\ot W),
\qquad u\ot v\ot w\longmapsto\omega\cdot(u\ot v\ot w),
\]
is a finite tensor category, denoted by $\Coh(G,\omega)$ \cite{G1}. The underlying tensor product is over $k$, with diagonal $\O(G)$-action. Its unit is $\delta_1$, not $\O(G)$.

For a normalized $2$-cochain $J$ on $G$, our convention is
\begin{equation}\label{dpsiiswomegaeq}
 dJ(a,b,c):=\frac{J(a,bc)J(b,c)}{J(ab,c)J(a,b)}.
\end{equation}
The identity underlying functor
\[
\mathsf I_J:\Coh(G,\omega\,dJ)\longrightarrow\Coh(G,\omega)
\]
is a tensor equivalence with structure isomorphism
\[
(\mathsf I_J)_2(U,V):\mathsf I_J(U)\otimes\mathsf I_J(V)
\longrightarrow\mathsf I_J(U\otimes V),\qquad z\longmapsto J\cdot z.
\]
Its coherence equation is precisely \eqref{dpsiiswomegaeq}.

For $g\in G(k)$, define
\begin{gather}\label{omega123}
\omega_1:=(g\ot\id\ot\id)(\omega),\qquad
\omega_2:=(\id\ot g\ot\id)(\omega),\\
\omega_3^g:=(\Ad_g\ot\Ad_g)(\id\ot\id\ot g)(\omega).
\end{gather}
Set
\begin{equation}\label{omegag}
\omega_g:=\omega_1(\Ad_g\ot\id)(\omega_2^{-1})\omega_3^g.
\end{equation}
Thus,
\[
\omega_g(a,b)=
\frac{\omega(g,a,b)\omega(gag^{-1},gbg^{-1},g)}{\omega(gag^{-1},g,b)}.
\]
This is a $2$-cochain on $G$, and its restriction to the centralizer of $g$ is a $2$-cocycle; see Lemma \ref{omegagprop}.

\subsection{(Bi)equivariant sheaves and (bi)comodules}\label{sec:Equivariant sheaves and comodules}
Let $H,K\subseteq G$ and $\psi,\eta$ be as in Definition \ref{defequivbi}. Regard $\O(H)$ and $\O(K)$ as sheaves on $G$ by extension by zero. Let
\[
\Delta_\psi(f):=\Delta(f)\psi,\qquad
\Delta_\eta(f):=\Delta(f)\eta.
\]

\begin{lemma}\label{internal-coalgebra}
The objects $\O(H)_\psi$ and $\O(K)_\eta$, with counits induced by the identity, are coalgebras in $\Coh(G,\omega)$. In particular,
\[
(\id\ot\Delta_\psi)\Delta_\psi
=\alpha_{\O(H),\O(H),\O(H)}
(\Delta_\psi\ot\id)\Delta_\psi.
\]
\end{lemma}

\begin{proof}
This follows since the two iterated comultiplications differ by
\[
(\id\ot\Delta)(\psi)(1\ot\psi)
\bigl((\Delta\ot\id)(\psi)(\psi\ot1)\bigr)^{-1}
=d\psi=\omega_{\mid H},
\]
and normalization gives the counit identities. 

The proof for $K$ is identical.
\end{proof}

\begin{definition}\label{defright}
An $(\O(H)_\psi,\O(K)_\eta)$-bicomodule in $\Coh(G,\omega)$ is an object $S$ with counital $\O(G)$-linear maps
\[
\lambda:S\to\O(H)_\psi\ot S,\qquad
\rho:S\to S\ot\O(K)_\eta
\]
satisfying
\begin{equation}\label{bicomodule-identities}
\begin{split}
(\id\ot\lambda)\lambda
&=\alpha_{H,H,S}(\Delta_\psi\ot\id)\lambda,\\
(\id\ot\Delta_\eta)\rho
&=\alpha_{S,K,K}(\rho\ot\id)\rho,\\
(\id\ot\rho)\lambda
&=\alpha_{H,S,K}(\lambda\ot\id)\rho.
\end{split}
\end{equation}
The subscripts $H,K$ on $\alpha$ abbreviate the corresponding sheaves. The left and right $\O(G)$-actions in these formulas are induced by $(h,x)\mapsto hx$, and $(x,k)\mapsto xk$. The one-sided comodule categories are defined by the corresponding single identity.
\end{definition}

\begin{proposition}\label{abelcat}
There is an abelian equivalence
\[
\Coh^{((H,\psi),(K,\eta))}(G,\omega)
\simeq {\rm Bicomod}_{\Coh(G,\omega)}(\O(H)_\psi,\O(K)_\eta),
\]
which preserves the underlying sheaf. The analogous statements hold for left and right equivariant sheaves.
\end{proposition}

\begin{proof}
The coactions give translation maps
\[
L_h(x):S_x\longrightarrow S_{h^{-1}x},\qquad
R_k(x):S_x\longrightarrow S_{xk^{-1}}.
\]
Counitality and \eqref{bicomodule-identities}, applied to an element and its inverse, make these maps invertible. The three identities become
\begin{align*}
L_{h'}(h^{-1}x)L_h(x)
 &=\psi(h,h')\omega(h,h',(hh')^{-1}x)L_{hh'}(x),\\
R_k(xk'^{-1})R_{k'}(x)
 &=\frac{\eta(k,k')}{\omega(x(kk')^{-1},k,k')}R_{kk'}(x),\\
R_k(h^{-1}x)L_h(x)
 &=\omega(h,h^{-1}xk^{-1},k)L_h(xk^{-1})R_k(x).
\end{align*}
Define
\begin{equation}\label{bicomodule-to-transport}
\theta_{(h,k)}(x):=\omega(h,h^{-1}x,k)
R_k(h^{-1}xk)^{-1}L_h(x):S_x\longrightarrow S_{h^{-1}xk}.
\end{equation}
For $a=(h,k)$, $b=(h',k')$ and $z=(hh')^{-1}x$, the unnormalized maps $R_k^{-1}L_h$ have composition multiplier
\[
\frac{\psi(h,h')}{\eta(k,k')}
\omega(h,h',z)\omega(h',z,k)\omega(z,k,k').
\]
Multiplying by the normalization factors in \eqref{bicomodule-to-transport} and using the pentagon gives \eqref{intro-tau}. Thus \eqref{transport-law} holds.

Conversely, recover the maps by
\[
L_h(x)=\theta_{(h,1)}(x),\qquad
R_k(x)=\theta_{(1,k)}(xk^{-1})^{-1}.
\]
The restrictions of \eqref{transport-law} to the two factors, and its mixed restriction, give the three displayed identities and hence \eqref{bicomodule-identities}. These constructions are inverse, and morphisms commute with $\theta$ exactly when they commute with both coactions. All formulas are identities of pulled-back sheaf maps.
\end{proof}

\subsection{Principal homogeneous spaces}\label{sec:Right principle homogeneous spaces}
Let $F$ act freely on a finite scheme $X$, and let $\pi:X\twoheadrightarrow X/F$ be the quotient morphism. Thus, $\pi$ is a faithfully flat $F$-torsor. Fix a normalized cocycle $\Phi$ as in \eqref{groupoid-cocycle}.

\begin{lemma}\label{helpful2-2}
The algebra $A(X,F,\Phi)$ (\ref{descent-algebra-definition}) is a split Azumaya algebra of rank $|F|^2$ over $\O(X/F)$. Thus, there is an invertible sheaf $E_\Phi$ on $X$, with $(F,\Phi)$-equivariant structure, such that
\[
\Coh(X/F)\xrightarrow{\simeq}\Coh^{(F,\Phi)}(X),\qquad U\longmapsto \pi^*U\ot_{\O(X)}E_\Phi
\]
is an equivalence. The choice of $E_\Phi$, including its equivariant structure, is part of this equivalence.
\end{lemma}

\begin{proof}
Let $A:=A(X,F,\Phi)$. Invariant functions embed $\O(X/F)$ centrally in $A$. After the faithfully flat base change $X\twoheadrightarrow X/F$, the torsor is trivial. Let $Y$ be its object scheme over the new base and $o$ its distinguished section. Write the pulled-back cocycle on the pair groupoid as $\varphi(y_0,y_1,y_2)$. The cocycle identity on $(o,y_0,y_1,y_2)$ gives
\[
\varphi(y_0,y_1,y_2)
=\frac{b(y_0,y_1)b(y_1,y_2)}{b(y_0,y_2)},\qquad
b(y_0,y_1):=\varphi(o,y_0,y_1).
\]
A gauge change therefore identifies the twisted algebra with the endomorphism algebra of $\O(Y)$ over the new base. Since $\O(Y)$ is free of rank $|F|$, $A$ is Azumaya of rank $|F|^2$ over $\O(X/F)$.

To prove splitting, consider a local factor $B$ of $\O(X/F)$, with maximal ideal $\mathfrak m$. The reduction of the corresponding algebra $A_B$ modulo $\mathfrak m$ is a matrix algebra over $k$. Lift a primitive matrix idempotent to $e\in A_B$; this is possible because $\mathfrak m A_B$ is nilpotent. The corner $eA_Be$ is a projective $B$-module whose reduction has dimension one, so $B\to eA_Be$ is an isomorphism by Nakayama's lemma. The idempotent is full, again by reduction and Nakayama, and Morita theory gives
\[
A_B\cong\End_B(A_Be).
\]
This splits $A$ on every local factor.

Regard the splitting module as an equivariant sheaf by Proposition~\ref{descent-algebra}. It is projective over $\O(X)$ since it is a direct summand of $A$. After the same faithfully flat base change, a splitting module for the matrix algebra is the standard module tensored with an invertible module over the base. Its restriction to $\O(Y)$ is therefore invertible, and invertibility descends. This gives $E_\Phi$. Morita equivalence is precisely $U\mapsto\pi^*U\ot_{\O(X)}E_\Phi$.
\end{proof}

Take now $X:=G$, $F:=H\subseteq G$, acting on the right, and $d\psi=\omega_{\mid H}$. The quotient is
\begin{equation}\label{freerightactionhong}
\pi:G\longrightarrow G/H,\qquad \mu:G\times H\to G,\quad (x,h)\mapsto xh
\end{equation}
where $\mu$ is the action, and its multiplier is
\[
\Phi(x;h,h'):=\psi(h,h')^{-1}\omega(x,h,h').
\]
Choose a splitting sheaf $E_\psi$ as in Lemma \ref{helpful2-2}; its underlying sheaf is isomorphic to $\O(G)$.

\begin{theorem}\label{helpful2}
There is an equivalence of abelian categories
\[
\Coh(G/H)\xrightarrow{\simeq}\Coh^{(H,\psi)}(G,\omega),
\qquad U\longmapsto\pi^*U\ot_{\O(G)}E_\psi.
\]
In particular, $\Coh^{(H,\psi)}(G,\omega)\simeq\Coh^{(H,1)}(G)$ as abelian categories.
\end{theorem}

\begin{proof}
Apply Lemma \ref{helpful2-2}. Tensoring an equivariant sheaf with $E_\psi^{-1}$ cancels its multiplier, and gives ordinary descent along $\pi$, hence it describes an inverse.
\end{proof}

\begin{remark}
The abelian equivalence in Theorem~\ref{helpful2} depends on a choice of $E_\psi$. Its underlying invertible sheaf is free, but its equivariance retains the twisting. The theorem does not assert that ordinary pullback is tensor for convolution. \qed
\end{remark}

For $\bar g:=gH\in(G/H)(k)$, set
\begin{equation}\label{varioussimplespre}
S_{\bar g}:=\pi^*\delta_{\bar g}\ot_{\O(G)} E_\psi,
\qquad P_{\bar g}:=P_{\Coh(G/H)}(\delta_{\bar g}).
\end{equation}

\begin{corollary}\label{simprojhpsi}
The objects $S_{\bar g}$, for $\bar g\in(G/H)(k)$, are all the simple objects of $\Coh^{(H,\psi)}(G,\omega)$, without repetitions. Moreover,
\[
P(S_{\bar g})=\pi^*P_{\bar g}\ot_{\O(G)} E_\psi.
\]
As sheaves on $G$, these objects are respectively $\O(gH)$, and $\bigoplus_{h\in H(k)}P_{gh}$. After choosing $S_{\bar1}$, and compatible identifications, convolution by $\delta_g$ gives
\[
S_{\bar g}\cong\delta_g\ot S_{\bar1},\qquad
P(S_{\bar g})\cong\delta_g\ot P(S_{\bar1}).
\]
\end{corollary}

\begin{proof}
An equivalence preserves simple objects and their projective covers. Pullback of the local factor of $\O(G/H)$ at $\bar g$ is the sum of the local factors of $\O(G)$ at $gh$, $h\in H(k)$. Tensoring by the invertible sheaf does not change these underlying modules up to isomorphism. Finally, $\delta_g$ is invertible in $\Coh(G,\omega)$, and its action on the module category translates the support from $H$ to $gH$. It therefore sends the unique simple on the first coset to that on the second, and sends its projective cover to the corresponding cover.
\end{proof}

The same argument applies to the free left action of $H$, and the quotient
\begin{equation}\label{freerightactionhong2}
p:G\longrightarrow H\backslash G.
\end{equation}

\subsection{\texorpdfstring{Twisting by the inverse $2$-cocycle}{Twisting by the inverse 2-cocycle}}\label{sec:Twisting by the inverse $2$-cocycle}
Let $A$ be a finite group scheme, and $\xi\in Z^2(A,\Gm)$. A left $\O(A)_\xi$-comodule gives operators $\ell_a$ satisfying
\[
\ell_b\ell_a=\xi(a,b)\ell_{ab},\qquad a,b\in A.
\]
The inverse operators $r_a:=\ell_a^{-1}$ satisfy
\[
r_ar_b=\xi(a,b)^{-1}r_{ab},\qquad a,b\in A.
\]
Thus, inversion of transport gives an equivalence
\begin{equation}\label{fromlefttoright}
{\rm Rep}(A,\xi)_k\simeq {\rm Rep}_k(A,\xi^{-1}),
\end{equation}
where the category on the right consists of right $\O(A)_{\xi^{-1}}$-comodules. These formulas hold on points over arbitrary commutative algebras, and hence give maps of comodules. More explicitly, if $\ell(v)=\sum v^{-1}\ot v^0$, and
\[
Q_\xi:=m(S\ot\id)(\xi),
\]
then the right coaction is
\[
v\longmapsto\sum v^0\ot S(v^{-1})Q_\xi^{-1}.
\]
Indeed, $Q_\xi(a)=\xi(a^{-1},a)$, and $\xi(a^{-1},a)=\xi(a,a^{-1})$ by the cocycle identity, so the displayed map is exactly inversion of $\ell_a$. In particular, a left comodule, and a right comodule with the \emph{same} multiplier, should not be identified without this inversion.

\subsection{\texorpdfstring{Module categories over $\Coh(G,\omega)$}{Module categories over Coh(G, omega)}}\label{sec:M(H,psi)}
Fix $G,\omega,H,\psi$ as above, and set
\begin{equation}\label{eq:M(H,psi)}
\M(H,\psi):=\Coh^{(H,\psi)}(G,\omega).
\end{equation}
It has a canonical structure of an indecomposable exact left module category over $\D:=\Coh(G,\omega)$, given by convolution \cite{G1}. In internal-comodule notation the action is
\[
X\ot^{\M}(S,\rho)=
\bigl(X\ot S,\,\alpha_{X,S,\O(H)}^{-1}(\id_X\ot\rho)\bigr).
\]
In particular, the associator is part of this coaction. Exactness and the classification below are the module-category results of \cite{EO,G1}; Theorem \ref{helpful2} describes only the underlying abelian category.

\begin{example}\label{simpexs0}
We have $\M(1,1)=\Coh(G,\omega)$, the regular module category. For $\omega=1$, $\M(G,1)\simeq{\rm Vec}$ is the usual fiber-functor module category of $\Coh(G)$. \qed
\end{example}

For a cochain $f$ on $G$, write
\[
f^g(a_1,\ldots,a_n):=f(ga_1g^{-1},\ldots,ga_ng^{-1}).
\]
The pentagon identities at
\[
\begin{gathered}
(g,a,b,c),\quad(gag^{-1},g,b,c),\\
(gag^{-1},gbg^{-1},g,c),\quad
(gag^{-1},gbg^{-1},gcg^{-1},g)
\end{gathered}
\]
give $d\omega_g=\omega/\omega^g$. Thus,
\begin{equation}\label{psig}
\psi_g:=(\psi^g\omega_g)|_{g^{-1}Hg}
\quad\text{satisfies}\quad d\psi_g=\omega|_{g^{-1}Hg}.
\end{equation}
Here $\psi^g(a,b)=\psi(gag^{-1},gbg^{-1})$ on $g^{-1}Hg$.

\begin{definition}\label{defequivpairs}
Two pairs $(H,\psi)$ and $(H',\psi')$ are said to be equivalent if $H'=g^{-1}Hg$ for some $g\in G(k)$, and $(\psi')^{-1}\psi_g$ is a coboundary on $H'$.
\end{definition}

The scalar correction in \eqref{psig} is essential; compare \cite[Theorem~1.1]{N} in the finite-group fusion setting.

\begin{theorem}\cite[Theorem~5.3]{G1}\label{modomega}
The assignment $(H,\psi)\mapsto\M(H,\psi)$ is a bijection between equivalence classes of pairs in Definition \ref{defequivpairs}, and equivalence classes of indecomposable exact left $\Coh(G,\omega)$-module categories. The equivalence with internal right $\O(H)_\psi$-comodules is a module equivalence.
\end{theorem}

\begin{remark}\label{helpful}
If $G$ is connected, then conjugation by closed points disappears, but equivalence modulo $2$-coboundaries does not. Thus, the classification still uses pairs modulo coboundaries, not individual cochains. \qed
\end{remark}

\begin{remark}\label{rem:N(H,psi)}
Similarly $\N(H,\psi):={}^{(H,\psi)}\Coh(G,\omega)$ is an indecomposable exact right module category. The left-comodule version of the preceding arguments gives the analogous classification. \qed
\end{remark}

\subsection{\texorpdfstring{Module categories over $\C(G,\omega,H,\psi)$}{Module categories over C(G, omega, H, psi)}}\label{sec:GCTC}
Let $\D:=\Coh(G,\omega)$ and $\M:=\M(H,\psi)$. We use
\begin{equation}\label{gsctcdef}
\C(G,\omega,H,\psi)=\Fun_\D(\M,\M)^{\rm rev}.
\end{equation}
Here $\Fun_\D$ denotes additive $k$-linear module functors and module natural transformations, and ${\rm rev}$ means opposite composition, not the opposite abelian category. Thus, $F\ot G=G\circ F$. This fixes consistently both the example $\C(G,\omega,1,1)\simeq\D$, and the left action by precomposition on $\Fun_\D(\M,\N)$. With the convention of ordinary composition all the corresponding tensor categories and module actions must instead be reversed.

\begin{example}\label{simpexs}
There are tensor equivalences
\[
\C(G,\omega,1,1)\simeq\Coh(G,\omega),\qquad
\C(G,1,G,1)\simeq{\rm Rep}_k(G).
\]
The symmetry of ${\rm Rep}_k(G)$ identifies it with its reverse. \qed
\end{example}

\begin{theorem}\cite[Theorem~7.12.16]{EGNO}\label{modgth}
The assignment
\[
\N\longmapsto\Fun_\D(\M(H,\psi),\N)
\]
gives an equivalence between the $2$-categories of indecomposable exact left $\D$-module categories, and of indecomposable exact left $\C(G,\omega,H,\psi)$-module categories.
\end{theorem}

For $(K,\eta)$ as above, set
\begin{equation}\label{mhpsieta}
\M((H,\psi),(K,\eta)):=\Coh^{((H,\psi),(K,\eta))}(G,\omega).
\end{equation}

\begin{lemma}\label{newsimple}
There are equivalences of abelian categories
\[
\begin{split}
\M((H,\psi),(K,\eta))
&\simeq {\rm Bicomod}_\D(\O(H)_\psi,\O(K)_\eta)\\
&\simeq\Fun_\D(\M(H,\psi),\M(K,\eta)).
\end{split}
\]
Under the second equivalence, a bicomodule $S$ acts by
\[
T\longmapsto T\mathbin{\square}_{\O(H)_\psi} S,
\]
where the cotensor product $\mathbin{\square}_{\O(H)_\psi}$ is the equalizer defined by the two coactions, with the ambient associator inserted to give a common target.
\end{lemma}

\begin{proof}
The first equivalence is Proposition~\ref{abelcat}. For the second, use the coalgebra form of the module-functor theorem; the algebra form is \cite[Proposition~7.11.1]{EGNO}, and rigidity gives the coalgebra version. A bicomodule defines the displayed cotensor functor. Conversely, a module functor $F$ gives $F(\O(H)_\psi)$, with left coaction obtained from the regular coaction using its module constraint, and right coaction inherited from $\M(K,\eta)$.

The cotensor functor is additive and is a $\D$-module functor because tensoring in $\D$ is exact and preserves the defining equalizer. It is therefore exact by \cite[Proposition~7.6.9]{EGNO}, since its source is an exact module category. Applying a module functor to the regular-comodule equalizer presentation identifies the two constructions as inverse module functors.
\end{proof}

\begin{lemma}\label{equdefofgstc}
The category $\C(G,\omega,H,\psi)$ is tensor equivalent to the category of $(\O(H)_\psi,\O(H)_\psi)$-bicomodules, with cotensor product. In particular, it is a finite tensor category.
\end{lemma}

\begin{proof}
If $F_S(T):=T\square S$, then $F_{S\square S'}=F_{S'}\circ F_S$, with the canonical associativity identifications. Thus, opposite composition in \eqref{gsctcdef} corresponds to the indicated cotensor product. Finiteness and rigidity are the dual-category theorem for exact module categories \cite{EO}.
\end{proof}

\section{Associated biequivariant sheaves}\label{S:Biequivariant sheaves on group schemes}

Let $B,C$ be finite group schemes and write $F:=B\times C$. Consider the canonical embeddings
\begin{equation}\label{injective morphisms}
j_1:B\to F,\quad b\mapsto(b,1),\qquad
j_2:C\to F,\quad c\mapsto(1,c).
\end{equation}
Let $\partial:A\hookrightarrow F$ be a closed subgroup scheme embedding. We can in particular take $\partial:=(\partial_1,\partial_2)$ for embeddings in each factor, as in the double-coset application. Let
\begin{equation}\label{embeddingdelta}
p:F\twoheadrightarrow A\backslash F,\qquad x_0:=p(1).
\end{equation}
The group $F$ acts on $A\backslash F$ by right multiplication. Fix a normalized action-groupoid cocycle $\Phi$ on $A\backslash F\times F\times F$. Set
\begin{equation}\label{zetafather}
\xi(a,b):=\Phi(x_0;\partial(a),\partial(b)),\qquad a,b\in A.
\end{equation}
Since $A$ fixes $x_0$, \eqref{groupoid-cocycle} implies $\xi\in Z^2(A,\Gm)$.
Let
\begin{equation}\label{CohBpsiCetaZomega}
\A:=
\Coh^{(F,\Phi)}(A\backslash F)=
\Coh^{(B\times C,\Phi)}(A\backslash (B\times C)).
\end{equation}

Consider the action of $A\times F$ on $F$ from the right by
\[
q\cdot(a,b)=\partial(a)^{-1}qb,\qquad a\in A,\quad q,b\in F,
\]
and define
\[
\widetilde\Phi\bigl(q;(a,b),(a',b')\bigr)
:=
\Phi(p(q);b,b'),\qquad a,a'\in A,\quad q,b,b'\in F.
\]
Since $p(q\cdot(a,b))=p(q)\cdot b$, the cocycle identity for $\Phi$ implies that for $\widetilde\Phi$. Let
\begin{equation}\label{CohBpsiCetaZomegaB}
\mathscr B:=
\Coh^{(A\times F,\widetilde\Phi)}(F)=\Coh^{(A\times B\times C,\widetilde\Phi)}(B\times C).
\end{equation}

Let $\lambda:A\times F\to F$ and $m:F\times F\to F$ be given by $\lambda(a,q)=\partial(a)q$ and $m(q,b)=qb$. Equivalently, an object of $\B$ is a triple $(E,D,T)$, where $E\in\Coh(F)$ and
\[
D:\operatorname{pr}_2^*E\xrightarrow{\cong}\lambda^*E,
\qquad T:\operatorname{pr}_1^*E\xrightarrow{\cong}m^*E.
\]
For a test point $q:\Spec(R)\to F$, write $E_q:=q^*E$. Suppressing $\partial$, the pullback maps $D_a(q):E_q\to E_{aq}$ and $T_b(q):E_q\to E_{qb}$ satisfy
\[
\begin{gathered}
D_1(q)=T_1(q)=\id_{E_q},\qquad
D_{a'}(aq)D_a(q)=D_{a'a}(q),\\
T_c(qb)T_b(q)=\Phi(p(q);b,c)T_{bc}(q),
\end{gathered}
\]
and
\begin{equation}\label{descent-compatibility}
D_a(qb)T_b(q)=T_b(aq)D_a(q).
\end{equation}
Thus $D$ is ordinary left $A$-descent, and $T$ is twisted right $F$-transport. Equation~\eqref{descent-compatibility} equates the two maps $E_q\to E_{aqb}$; it is exactly the condition that $T$ descend along $p\times\id_F$.

The transport for \eqref{CohBpsiCetaZomegaB} is
\[
U_{(a,b)}(q)=D_{a^{-1}}(qb)T_b(q).
\]
The preceding identities give its composition law with multiplier $\widetilde\Phi$. Conversely, restriction gives $D_a(q)=U_{(a^{-1},1)}(q)$ and $T_b(q)=U_{(1,b)}(q)$. Morphisms are sheaf maps commuting with both $D$ and $T$.

\subsection{\texorpdfstring{The equivalence $\B\simeq{\rm Rep}(A,\xi)_k$}{The equivalence B and Rep(A, xi)}}
For $q,b\in F$, set 
\[\Phi_0(q,b):=\Phi(x_0;q,b).\]
Let $(V,\ell)\in{\rm Rep}(A,\xi)_k$, so that $\ell_b\ell_a=\xi(a,b)\ell_{ab}$. On the free sheaf 
\[E_V:=\O(F)\ot_kV\]
define
\begin{equation}\label{explicit-descent}
D_a(q):=\Phi_0(a,q)\ell_a^{-1},\qquad
T_b(q):=\Phi_0(q,b)\id_V.
\end{equation}
The inclusion $\partial$ is suppressed in these formulas.

\begin{lemma}\label{lemmatezer}
The maps in \eqref{explicit-descent} make $E_V$ an object of $\B$.
\end{lemma}

\begin{proof}
The left descent identity is $D_b(aq)D_a(q)=D_{ba}(q)$. Its scalar part is
\[
\Phi_0(b,aq)\Phi_0(a,q)=\xi(b,a)\Phi_0(ba,q),
\]
which is \eqref{groupoid-cocycle} at $(x_0;b,a,q)$, since $A$ fixes $x_0=p(1)$. Its operator part is
$\ell_b^{-1}\ell_a^{-1}=\xi(b,a)^{-1}\ell_{ba}^{-1}$.
The right transport identity follows from
\[
\Phi_0(q,b)\Phi_0(qb,c)=\Phi(p(q);b,c)\Phi_0(q,bc).
\]
Finally, the compatibility \eqref{descent-compatibility} is
\[
\Phi_0(a,qb)\Phi_0(q,b)=\Phi_0(a,q)\Phi_0(aq,b),
\]
again the cocycle identity, this time at $(x_0;a,q,b)$. Normalization gives the identity maps.
\end{proof}

\begin{theorem}\label{lemmafromg}
The assignment
\[
\Rep(A,\xi)_k\longrightarrow\B,\qquad
V\longmapsto E_V=\O(F)\ot_kV,
\]
with structure \eqref{explicit-descent}, is an equivalence.
\end{theorem}

\begin{proof}
Given $E\in\B$, let $V:=E_1$. Right transport identifies $E$ with the free sheaf $\O(F)\ot V$ using $T_q(1):V\to E_q$. Define
\[
\ell_a:=D_a(1)^{-1}T_a(1):V\longrightarrow V.
\]
Then the commuting descent and right-transport identities show that $\ell_b\ell_a=\xi(a,b)\ell_{ab}$. In the trivialization by $T_q(1)$, the right transport is $\Phi_0(q,b)\id$ by its composition law. The left descent is $\Phi_0(a,q)\ell_a^{-1}$, by \eqref{descent-compatibility}. Thus, it is exactly \eqref{explicit-descent}. A morphism commuting with right transport is determined by its value at $1$, and commuting with left descent is precisely the condition that this value commute with $\ell$. This proves full faithfulness and essential surjectivity.
\end{proof}

\subsection{\texorpdfstring{The equivalence $\A\simeq\B$}{The equivalence A and B}}
The quotient $p:F\to A\backslash F$ is a faithfully flat left $A$-torsor. Pullback gives an ordinary descent datum on $p^*S$ for any sheaf $S$ on $A\backslash F$. A right $F$-transport on $S$ pulls back to a right transport of multiplier $p^*\Phi$ compatible with this datum.

\begin{theorem}\label{flat-descent-equivalence}
Pullback and faithfully flat descent give inverse equivalences
\[
p^*:\A\xrightarrow{\simeq}\B,\qquad
p_*^A:\B\xrightarrow{\simeq}\A.
\]
\end{theorem}

\begin{proof}
Faithfully flat descent \cite[Proposition~35.5.2]{Stacks} identifies coherent sheaves on $X$ with sheaves on $F$ carrying the stated left $A$-descent datum. Apply the same descent to the right transport map. The condition for it to descend is exactly \eqref{descent-compatibility}. The identity \eqref{transport-law} descends because it can be checked after the faithfully flat pullback. Morphisms and their compositions descend in the same way.
\end{proof}

\subsection{\texorpdfstring{The first equivalence ${\rm Rep}(A,\xi)_k\simeq\A$}{The first equivalence Rep(A, xi) and A}}\label{sec:The first equivalence}
Set
\begin{equation}\label{indcomod}
{\rm Ind}(V):=p_*^A(E_V),\qquad
{\rm Res}(S):=S_{x_0},
\end{equation}
where the fiber has its restricted stabilizer transport.

\begin{theorem}\label{modrep0-1}
There are mutually inverse equivalences of abelian categories
\[
{\rm Ind}:{\rm Rep}(A,\xi)_k\xrightarrow{\simeq}\A=\Coh^{(F,\Phi)}(A\backslash F),
\qquad {\rm Res}(S)=S_{x_0}.
\]
For each $V$, $p^*{\rm Ind}(V)$ is the free sheaf $\O(F)\ot V$ with the maps \eqref{explicit-descent}. In particular, ${\rm Ind}(V)$ is locally free of rank $\dim(V)$ on $A\backslash F$.
\end{theorem}

\begin{proof}
Combine Theorems~\ref{lemmafromg} and \ref{flat-descent-equivalence}. Local freeness and its rank descend from $\O(F)\ot_k V$ along the faithfully flat morphism $p$.
\end{proof}

\subsection{\texorpdfstring{The second equivalence ${\rm Rep}(A,\xi)_k\simeq\A$}{The second equivalence Rep(A, xi) and A}}\label{sec:thesecondequivalence}
Every finite locally free module of constant rank over the finite-dimensional semilocal algebra $\O(A\backslash F)$ is free. For each $V$, choose an $\O(A\backslash F)$-linear isomorphism
\[
t_V:{\rm Ind}(V)\xrightarrow{\cong}\O(A\backslash F)\ot_k V.
\]
Transport the equivariant structure by $t_V$, and denote the resulting object by $\mathbf F(V)$. For a morphism $f:V\to W$, define
\[
\mathbf F(f):=t_W\,{\rm Ind}(f)\,t_V^{-1}.
\]

\begin{theorem}\label{modrep00}
The above choices give an equivalence of abelian categories
\[
\mathbf F:{\rm Rep}(A,\xi)_k\xrightarrow{\simeq}\Coh^{(F,\Phi)}(A\backslash F),
\qquad V\longmapsto\O(A\backslash F)\ot_k V.
\]
Its equivariance maps and its maps on morphisms are the transported maps just specified. In particular, the displayed formula is a formula for the underlying sheaf, not a canonical untwisted tensor-product functor.
\end{theorem}

\begin{proof}
The maps $t_V$ give a natural isomorphism from ${\rm Ind}$ to the functor $\mathbf F$ just defined, including its definition on morphisms. Theorem \ref{modrep0-1} therefore applies.
\end{proof}

\begin{remark}
The choices $t_V$ need not make $\mathbf F(f)=\id\ot f$. For computations, use \eqref{explicit-descent} and then choose bases. \qed
\end{remark}

\section{\texorpdfstring{Double cosets in $G$ and biequivariant sheaves on $(G,\omega)$}{Double cosets in G and biequivariant sheaves on (G, omega)}}\label{sec:doublecosetsandbiequivariant}
Fix a finite group scheme $G$, and $\omega\in Z^3(G,\Gm)$.

\subsection{\texorpdfstring{Double cosets in $G$}{Double cosets in G}}\label{sec:doublecosets}
For the reader's convenience, we recall the double-coset construction of \cite[Section 4]{GS}. Let $\iota_H:H\hookrightarrow G$ and $\iota_K:K\hookrightarrow G$ be closed subgroup scheme embeddings. Consider the right action
\begin{equation}\label{hkaction}
\mu_{G\times H\times K}:G\times(H\times K)\to G,
\qquad (x,h,k)\longmapsto h^{-1}xk.
\end{equation}
On coordinate algebras it is
\begin{equation}\label{hkactioncom}
\mu^{\sharp}(f)=f_2\ot\iota_H^{\sharp}S(f_1)\ot\iota_K^{\sharp}(f_3).
\end{equation}
Let $Y:=G/(H\times K)$ be the finite quotient scheme. Its closed points correspond to the orbits $Z_g:=HgK$, $g\in G(k)$. Thus,
\begin{equation}\label{doublecosetfadec}
\Coh(Y)=\bigoplus_{Z\in Y(k)}\Coh(Y)_Z.
\end{equation}
We distinguish the point $Z\in Y(k)$ from its orbit subscheme $Z\subseteq G$ when taking coordinate algebras or ideals.

Fix $g\in G(k)$, and let $L^g:=H\cap gKg^{-1}$, with intersection taken as a fiber product of group schemes. The embedding
\begin{equation}\label{thetag}
\partial_g:L^g\hookrightarrow H\times K,
\qquad l\longmapsto(l,g^{-1}lg)
\end{equation}
is the stabilizer of $g$. Its free left action on $H\times K$ is
\begin{equation}\label{lgactsonhk}
(l,h,k)\longmapsto(lh,g^{-1}lgk).
\end{equation}
There is an isomorphism of schemes, equivariant for the right $H\times K$-actions,
\begin{equation}\label{doucosetsch}
j_g:L^g\backslash(H\times K)\xrightarrow{\cong} Z_g,
\qquad [(h,k)]\longmapsto h^{-1}gk.
\end{equation}
In particular,
\[
|Z_g|=\frac{|H||K|}{|L^g|},\qquad
|Z_g(k)|=\frac{|H(k)||K(k)|}{|L^g(k)|}.
\]
The isomorphism follows by identifying the fibers of the orbit map with the left $L^g$-torsor fibers. It is a statement about the homogeneous scheme, not only its closed points.

\subsection{\texorpdfstring{Biequivariant sheaves on $(G,\omega)$}{Biequivariant sheaves on (G, omega)}}\label{sec:BiequivariantsheavesonGomega}
Suppose $\psi\in C^2(H,\Gm)$ and $\eta\in C^2(K,\Gm)$ satisfy $d\psi=\omega|_H$ and $d\eta=\omega|_K$. For $a=(h,k)$, $b=(h',k')$ and $x_1=x\cdot a$, $x_2=x\cdot(ab)$, define the multiplier on \eqref{hkaction} by
\begin{equation}\label{intro-tau}
\Phi(x;a,b):=
\frac{\psi(h,h')}{\eta(k,k')}
\frac{\omega(x,k,k')\omega(h,h',x_2)}{\omega(h,x_1,k')}.
\end{equation}
All factors are pulled back by the indicated scheme morphisms.

\begin{proposition}\label{taulemma}
The function \eqref{intro-tau} is a normalized action-groupoid $2$-cocycle. Its restriction at $g$ to $\partial_g(L^g)$ is
\begin{equation}\label{defnrmwg}
\begin{split}
\xi_g(l,m)&:=\psi(l,m)\eta(g^{-1}lg,g^{-1}mg)^{-1}w_g(l,m),\\
w_g(l,m)&:=\frac{\omega(g,g^{-1}lg,g^{-1}mg)\omega(l,m,g)}
{\omega(l,g,g^{-1}mg)}.
\end{split}
\end{equation}
In particular, $\xi_g\in Z^2(L^g,\Gm)$.
\end{proposition}

\begin{proof}
Normalization follows from normalization of $\omega,\psi,\eta$. We verify the full composition identity. Write $a_i:=(h_i,k_i)$, and $x_i:=h_i^{-1}x_{i-1}k_i$, with $x_0:=x$. Thus, $x_{i-1}k_i=h_i x_i$. Let
\[
W(x;a,b):=\frac{\omega(x,k,k')\omega(h,h',x\cdot ab)}
{\omega(h,x\cdot a,k')}.
\]
Expanding its four factors gives
\begin{equation}\label{omega-factor-ratio}
\frac{W(x;a_1,a_2)W(x;a_1a_2,a_3)}
{W(x_1;a_2,a_3)W(x;a_1,a_2a_3)}
=\frac{\omega(h_1,h_2,h_3)}{\omega(k_1,k_2,k_3)}.
\end{equation}
For completeness, the cancellations in \eqref{omega-factor-ratio} use the $3$-cocycle identity successively on
\[
(x,k_1,k_2,k_3),\quad(h_1,x_1,k_2,k_3),\quad
(h_1,h_2,x_2,k_3),\quad(h_1,h_2,h_3,x_3).
\]
At each step replace $x_{i-1}k_i$ by $h_i x_i$. Then all factors involving some $x_i$ cancel, and the two remaining factors are the right side of \eqref{omega-factor-ratio}. The analogous ratio for $\psi(h,h')/\eta(k,k')$ is
\[
(d\psi)(h_1,h_2,h_3)^{-1}(d\eta)(k_1,k_2,k_3),
\]
so their product is $1$. This is \eqref{groupoid-cocycle}. On the stabilizer, both intermediate points are $g$, and restriction gives \eqref{defnrmwg}. 

The groupoid identity then becomes the ordinary cocycle identity for $\xi_g$.
\end{proof}

\begin{proposition}\label{gauge-covariance}
Let $J\in C^2(G,\Gm)$ be normalized, and set
\[
\omega':=\omega dJ,\quad \psi':=\psi J_{\mid H},\quad
\eta':=\eta J_{\mid K}.
\]
For $a:=(h,k)$, define the normalized groupoid $1$-cochain
\[
b(x;a):=\frac{J(h,x\cdot a)}{J(x,k)}.
\]
Then
\begin{equation}\label{full-gauge}
\Phi'(x;a,b')=\Phi(x;a,b')
\frac{b(x;a)b(x\cdot a;b')}{b(x;ab')}.
\end{equation}
In particular, for $l\in L^g$ let
\begin{equation}\label{stabilizer-gauge}
b_g(l):=\frac{J(l,g)}{J(g,g^{-1}lg)}.
\end{equation}
Then $\xi'_g=\xi_g\,db_g$. Multiplication of transport by $b$ gives an equivalence of the corresponding equivariant categories, preserving underlying sheaves and projective covers.
\end{proposition}

\begin{proof}
Substitute $\omega dJ,\psi J_H,\eta J_K$ into \eqref{intro-tau}. In the resulting quotient, expand
\[
dJ(u,v,w)=\frac{J(u,vw)J(v,w)}{J(uv,w)J(u,v)}.
\]
With $x_1:=h^{-1}xk$, and $x_2:=(hh')^{-1}xkk'$, use $hx_1=xk$ and $h'x_2=x_1k'$. Then all factors except
\[
\frac{J(h,x_1)J(h',x_2)J(x,kk')}
{J(x,k)J(x_1,k')J(hh',x_2)}
\]
cancel. This is the quotient in \eqref{full-gauge}. Restriction gives \eqref{stabilizer-gauge}. Finally, $\theta'_a=b(x;a)\theta_a$ satisfies the primed transport identity by direct substitution. Since its inverse uses $b^{-1}$, the functor is an equivalence.
\end{proof}

\begin{proposition}\label{representative-change}
Let $r=(h_0,k_0)\in(H\times K)(k)$ and $g'=g\cdot r$. Set
\[
A_g:=\partial_g(L^g),\qquad A_{g'}:=\partial_{g'}(L^{g'}).
\]
Conjugation gives an isomorphism $\kappa:A_{g'}\to A_g$, $a'\mapsto ra'r^{-1}$. For
\[
u_r(a'):=\frac{\Phi(g;r,a')}{\Phi(g;\kappa(a'),r)},
\]
the restricted cocycles satisfy
\[
\xi_{g'}=\kappa^*\xi_g\,du_r,
\]
where the cocycles are identified with those on the embedded stabilizers. Under $A_g\cong L^g$ and $A_{g'}\cong L^{g'}$, the map $\kappa$ is $l'\mapsto h_0l'h_0^{-1}$. The change of stabilizer label is this conjugation followed by the indicated coboundary change.
\end{proposition}

\begin{proof}
For a sheaf $S$, compare the two factorizations of the arrow $ra'=\kappa(a')r$. By \eqref{transport-law},
\[
\theta_{a'}(g')=u_r(a')\theta_r(g)\theta_{\kappa(a')}(g)\theta_r(g)^{-1}.
\]
Composing this equality for two elements of $A_{g'}$ gives exactly the claimed identity of cocycles. The same identity follows by three applications of \eqref{groupoid-cocycle}, so it is an identity of functions independent of the choice of $S$.
\end{proof}

For $Z\in Y(k)$, choose $g\in Z(k)$, and use \eqref{doucosetsch} to transport Theorems \ref{modrep0-1} and \ref{modrep00}. Denote the resulting functors by
\begin{equation}\label{boldfzg1}
{\rm Ind}_Z,\ \mathbf F_Z:{\rm Rep}(L^g,\xi_g)_k\longrightarrow
\Coh^{((H,\psi),(K,\eta))}(Z).
\end{equation}
The target always means the restricted action-groupoid datum.

\begin{corollary}\label{newimpcor}
The functors in \eqref{boldfzg1} are equivalences. The first is the descent construction \eqref{explicit-descent}, and the second has underlying sheaf $\O(Z)\ot_kV$, with transported structure. Evaluation at $g$ is an inverse of the first equivalence.
\end{corollary}

\begin{proof}
The map $j_g$ is equivariant, and it identifies stabilizers and the complete multipliers. Apply Theorems \ref{modrep0-1} and \ref{modrep00}.
\end{proof}

\begin{remark}\label{orbit-not-augmented}
The restriction used is on $Z\times(H\times K)^2$, not on $Z^3$. \qed
\end{remark}

\section{\texorpdfstring{Module categories over $\C(G,\omega,H,\psi)$}
{Module categories over C(G, omega, H, psi)}}\label{sec:Module categories over GCTC}

Fix $\C:=\C(G,\omega,H,\psi)$, and an indecomposable exact left module category
\[
\M:=\M((H,\psi),(K,\eta))=\Coh^{((H,\psi),(K,\eta))}(G,\omega).
\]
Set $F:=H\times K$ and $Y:=G/F$, and let $\Phi$ be the multiplier \eqref{intro-tau}.

\subsection{\texorpdfstring{The simple objects of $\M((H,\psi),(K,\eta))$}{The simple objects of M((H, psi), (K, eta))}}\label{sec:The simple objects}
For $Z\in Y(k)$, let $\iota_Z:Z\hookrightarrow G$ be the orbit embedding and $\mathcal I(Z)\subseteq\O(G)$ its defining ideal. Set
\begin{equation}\label{defnofmzg}
\M_Z:=\{S\in\M\mid\mathcal I(Z)S=0\}.
\end{equation}
Then $\M_Z$ is the essential image of the fully faithful exact functor
\begin{equation}\label{mg}
\iota_{Z*}:\Coh^{((H,\psi),(K,\eta))}(Z)\longrightarrow\M.
\end{equation}
Let $\overline{\M_Z}$ be its Serre closure in $\M$, namely the full subcategory of objects all of whose composition factors lie in $\M_Z$.

For $g\in Z(k)$, compose the functors of Corollary \ref{newimpcor} with $\iota_{Z*}$. We use the same notation for the resulting functors into $\M$:
\begin{equation}\label{Indz}
{\rm Ind}_Z,\ \mathbf F_Z:{\rm Rep}(L^g,\xi_g)_k\longrightarrow\M.
\end{equation}

\begin{theorem}\label{simmodrep}
The following hold.
\begin{enumerate}
\item For every $Z\in Y(k)$, and $g\in Z(k)$, the functors ${\rm Ind}_Z$ and $\mathbf F_Z$ give equivalences ${\rm Rep}(L^g,\xi_g)_k\simeq\M_Z$.
\item The simple objects of $\M$ are ${\rm Ind}_Z(V)$, where $Z\in Y(k)$ and $V$ is a simple stabilizer representation. After fixing one $g$ in each orbit these labels have no repetitions.
\item There is a direct sum decomposition of abelian categories
\[
\M=\bigoplus_{Z\in Y(k)}\overline{\M_Z}.
\]
\end{enumerate}
\end{theorem}
\begin{proof}
Part (1) is Corollary \ref{newimpcor} and Proposition \ref{equivmorphs}. Here is a support argument proving both other parts, including the absence of cross extensions in the equivariant category.

Let $X_Z\subseteq G$ be the open-and-closed subscheme which is the union of the full local components at the points of $Z(k)$. Its idempotent $e_Z\in\O(G)$ is invariant: the equality $\mu^{\sharp}(e_Z)=e_Z\ot1$ can be checked on the closed points of $G\times F$, since idempotents of a finite algebra are uniquely determined modulo its nilradical. Therefore $e_ZS$ is an equivariant direct summand of every $S\in\M$, and
\[
S=\bigoplus_Ze_ZS.
\]
In $\O(X_Z)$, the ideal $I_Z$ defining $Z$ is invariant and nilpotent, since $Z$ and $X_Z$ have the same closed points. For a simple $S$ on $X_Z$, the subobject $I_ZS$ is either $0$ or $S$; nilpotence rules out the latter. Thus every such simple lies in $\M_Z$. Conversely every object on $X_Z$ has only these composition factors. This proves that the category on $X_Z$ is precisely $\overline{\M_Z}$ and gives the direct sum decomposition. Simplicity inside $\M_Z$ is the same as simplicity in $\M$, because every subobject of an object annihilated by $\mathcal I(Z)$ is also annihilated by that ideal.
\end{proof}

\begin{remark}
An orbit-indexed Serre summand need not be an indecomposable block. For example, on a point the category of representations of a finite group can have several blocks. We use ``Serre summand'' rather than ``block'' unless indecomposability has been proved. \qed
\end{remark}

\subsection{\texorpdfstring{Projectives in $\M((H,\psi),(K,\eta))$}{Projectives in M((H, psi), (K, eta))}}\label{sec:The projective objects}
Fix $Z$ and $g$ as above. Write $R_Z:=\O(X_Z)$, and let $I_Z\subseteq R_Z$ be the nilpotent ideal defining $Z$. Apply Proposition~\ref{descent-algebra} and Corollary~\ref{algebra-base-change} to obtain
\begin{equation}\label{orbit-algebras}
\begin{split}
A_Z&:=A(X_Z,H\times K,\Phi),\qquad J_Z:=I_ZA_Z,\\
\overline A_Z:=A_Z/J_Z&=A(Z,H\times K,\Phi_{\mid Z}).
\end{split}
\end{equation}
Thus, $A_Z\text{-}{\rm mod}\simeq\overline{\M_Z}$ and $\overline A_Z\text{-}{\rm mod}\simeq\M_Z$, with their underlying sheaves preserved. The ideal $J_Z$ is nilpotent, since it is generated by the invariant ideal $I_Z$.

\begin{proposition}\label{closedptproj}
Let $V$ be a simple object of ${\rm Rep}(L^g,\xi_g)_k$. Choose a primitive idempotent $\bar e$ of $\overline A_Z$, such that
\[
\overline A_Z\bar e\cong{\rm Ind}_Z\bigl(P_{(L^g,\xi_g)}(V)\bigr)
\]
as objects of $\M_Z$. Then $\bar e$ lifts to a primitive idempotent $e\in A_Z$, and
\[
A_Ze\cong P_{\M}({\rm Ind}_Z(V)).
\]
This object is independent, up to isomorphism, of the choices of the two idempotents.
\end{proposition}

\begin{proof}
The equivalence of Theorem~\ref{simmodrep}(1) preserves projective covers inside $\M_Z$, so such a primitive $\bar e$ exists. Idempotents lift modulo the nilpotent ideal $J_Z$. Explicitly, if $r=a^2-a$ is the defect of a lift, replacing $a$ by $a+(1-2a)r$ makes the new defect a multiple of $r^2$; iteration terminates. This works in every characteristic. A lift of a primitive idempotent is primitive, since a decomposition would either descend to a nontrivial decomposition of $\bar e$, or produce a nonzero idempotent in the nilpotent ideal $J_Z$.

The left ideal $A_Ze$ is projective. Every simple $A_Z$-module is annihilated by $J_Z$, so for any simple $W\in{\rm Rep}(L^g,\xi_g)_k$,
\[
\begin{split}
\Hom_{A_Z}(A_Ze,{\rm Ind}_Z(W))
&\cong\Hom_{\overline A_Z}(\overline A_Z\bar e,{\rm Ind}_Z(W))\\
&\cong\Hom_{{\rm Rep}(L^g,\xi_g)_k}(P(V),W).
\end{split}
\]
The dimension is $1$ for $W\cong V$, and $0$ otherwise. Morphisms to simple objects in a different Serre summand are zero by Theorem~\ref{simmodrep}. Thus the head is ${\rm Ind}_Z(V)$, proving the projective-cover assertion. Uniqueness gives independence of choices.
\end{proof}

\begin{theorem}\label{prsimmodrep}
For each $Z\in Y(k)$, representative $g\in Z(k)$, and simple $V\in{\rm Rep}(L^g,\xi_g)_k$, the following hold.
\begin{enumerate}
\item The projective cover is the object $A_Ze$ of Proposition \ref{closedptproj}. Its full equivariant structure is determined by \eqref{projective-coaction} below.
\item There is a (noncanonical) isomorphism of underlying $\O(G)$-modules
\begin{equation}\label{underlying-projective}
P_\M(\mathbf F_Z(V))\cong
\O(X_Z)\ot_k P_{(L^g,\xi_g)}(V).
\end{equation}
Since multiplication identifies
$G^{\circ}\times Z(k)$ with $X_Z$, the right side can also be written
\[
\O(G^{\circ})\ot\O(Z(k))\ot_kP_{(L^g,\xi_g)}(V),
\]
with $\O(G)$-action pulled back along that multiplication map.
\item In particular,
\[
\dim_kP_\M(\mathbf F_Z(V))
=|G^{\circ}|\,|Z(k)|\dim_kP_{(L^g,\xi_g)}(V).
\]
\end{enumerate}
\end{theorem}

\begin{proof}
Only the underlying-module assertion remains to be proved. By Proposition \ref{descent-algebra}, $A_Z$ is free as a left $R_Z$-module, therefore the summand $A_Ze$ is a finite projective $R_Z$-module. Its reduction modulo $I_Z$ is
\[
(A_Ze)/I_Z(A_Ze)\cong\overline A_Z\bar e
\cong{\rm Ind}_Z(P_{(L^g,\xi_g)}(V)).
\]
By Theorem \ref{modrep0-1}, this reduction is locally free over $\O(Z)$ of constant rank $r=\dim_kP_{(L^g,\xi_g)}(V)$. The subscheme $Z$ contains every closed point of $X_Z$. Thus, the rank of the projective $R_Z$-module $A_Ze$ is $r$ on every local factor. Since projective modules over each local Artin factor are free,  $A_Ze\cong R_Z^{\oplus r}$. This proves \eqref{underlying-projective}. Translation gives the asserted scheme isomorphism $G^{\circ}\times Z(k)\cong X_Z$, and then the dimension formula follows.
\end{proof}

The full equivariant structure is obtained as follows. Let 
$$\cR_Z:=\O(X_Z\times(H\times K))$$
be the coring used in Proposition \ref{descent-algebra}, with a right $R_Z$-basis $c_i$, and dual basis $f_i\in A_Z$. On $A_Ze$ its coaction is
\begin{equation}\label{projective-coaction}
ae\longmapsto\sum_i c_i\ot_{R_Z}(f_i\star a)e.
\end{equation}
The comultiplication of $\cR_Z$ uses the full function $\Phi$ on $X_Z$, not its value at a closed point. Formula \eqref{bicomodule-to-transport}, restricted to the two factors, recovers the left and right bicomodule maps. Choosing an $R_Z$-basis of $A_Ze$ transports these maps to \eqref{underlying-projective}.

\begin{remark}\label{free-sheaf-convention}
Every free-sheaf formula below specifies an underlying $\O(G)$-module. Its equivariance is obtained from $A_Z$ by \eqref{projective-coaction} and transported by the chosen isomorphism. No equivariant splitting or untwisted structure on a tensor factor is implicit. \qed
\end{remark}

\begin{example}\label{diagonal-splitting-warning}
Let $G:=\alpha_p$, and take $Z=G$. In coordinates the diagonal $\O(G)$-action on $k[x,y]/(x^p,y^p)$ is $t\mapsto x+y$. The map
\[
\chi:k[x,y]/(x^p,y^p)\to k[x]/(x^p),\qquad f(x,y)\mapsto f(x,0)
\]
is $\O(G)$-linear, but its algebra section $\nu(f(x))=f(x)$ is not since $t\nu(1)=x+y$, whereas $\nu(t\cdot1)=x$. An $\O(G)$-linear section in this example is $f(x)\mapsto f(x+y)$. Thus the section used in \cite[Lemma~5.4(1)]{GS} is not $\O(G)$-linear under the diagonal action. Proposition~\ref{closedptproj} replaces this splitting argument by a summand in the full equivariant algebra. \qed
\end{example}

\begin{example}\label{structcohhpsi}
If $K=1$, then $Y=H\backslash G$, $L^g=1$, and the unique simple on $Hg$ has underlying sheaf $\O(Hg)$. Its projective cover has underlying sheaf
\[
\O(X_{Hg})=\bigoplus_{h\in H(k)}P_{hg}.
\]
Its left equivariant structure is the descent structure of \eqref{projective-coaction}, or equivalently that of the left-hand version of Theorem \ref{helpful2}. \qed
\end{example}

\begin{corollary}\label{orbit-filtration}
Fix $Z$ and its representative $g$, and let $V\in\Rep(L^g,\xi_g)_k$ be simple. Set $R_Z:=\O(X_Z)$, let $I_Z$ be the ideal of $Z$ in $X_Z$, and let $P:=P_\M(\Ind_Z(V))$. For every $n\ge0$, multiplication induces an equivariant isomorphism
\[
I_Z^nP/I_Z^{n+1}P
\cong (I_Z^n/I_Z^{n+1})\otimes_{\O(Z)}
\Ind_Z\bigl(P_{(L^g,\xi_g)}(V)\bigr).
\]
The first factor has its ordinary equivariance, and the second has the restricted twisted equivariance.
\end{corollary}

\begin{proof}
The underlying $R_Z$-module of $P$ is projective by Theorem \ref{prsimmodrep}. Tensoring the exact sequence
\[0\to I_Z^{n+1}\to I_Z^n\to I_Z^n/I_Z^{n+1}\to0\]
with $P$ therefore gives
\[
I_Z^nP/I_Z^{n+1}P
\cong (I_Z^n/I_Z^{n+1})\otimes_{R_Z}P.
\]
The right-hand side depends only on $P/I_ZP$, which is the orbit-category projective cover in Proposition \ref{closedptproj}. Multiplication is equivariant because the transport on $P$ is semilinear over ordinary transport on $R_Z$. This gives the stated isomorphism, including its multiplier.
\end{proof}

This is the $I_Z$-adic filtration, not necessarily the radical filtration. For the algebra $B_\lambda$ of Theorem~\ref{alpha2center} below, the orbit ideal is $(X)$ and both cases have two $(X)$-adic layers of dimension two, although their radical layers differ as in Corollary~\ref{alpha2extensions}.

\subsection{\texorpdfstring{Fiber functors on $\C(G,\omega,H,\psi)$}{Fiber functors on C(G, omega, H, psi)}}
Recall that a fiber functor corresponds to an exact module category equivalent to ${\rm Vec}$.

\begin{corollary}\label{fibfungth}
Equivalence classes of fiber functors on $\C(G,\omega,H,\psi)$ correspond to equivalence classes of pairs $(K,\eta)$, such that
\[
d\eta=\omega_{\mid K},\qquad HK=G,
\]
and $\xi_1\in Z^2(H\cap K,\Gm)$ is nondegenerate. Here, $HK=G$ is equality of schemes, and nondegenerate means that the dual algebra of $\O(H\cap K)_{\xi_1}$ is central simple, equivalently its comodule category is ${\rm Vec}$.
\end{corollary}

\begin{proof}
If $HK=G$, there is just the orbit $G$, and Theorem~\ref{simmodrep} identifies the whole category with ${\rm Rep}(H\cap K,\xi_1)_k$. This is ${\rm Vec}$ precisely under the stated nondegeneracy condition.

Conversely suppose $\M((H,\psi),(K,\eta))\simeq{\rm Vec}$. Its unique simple belongs to one orbit category. Since every object is a direct sum of that simple, the defining ideal of this orbit annihilates every object. The regular module of $A(G,H\times K,\Phi)$ is faithful on $\O(G)$, by its embedded copy of $\O(G)$ in Proposition \ref{descent-algebra}. The orbit ideal is therefore zero, so the orbit is all of $G$, and since it contains $1$, it is $HK$. The same transitive equivalence now implies nondegeneracy. The equivalence relation on pairs and its compatibility with module equivalence are those of Theorems \ref{modomega} and \ref{modgth}.
\end{proof}

\section{\texorpdfstring{The structure of $\C(G,\omega,H,\psi)$}{The structure of C(G, omega, H, psi)}}\label{S:structure-GSC}

Fix $\C:=\C(G,\omega,H,\psi)$, and let $H^g:=H\cap gHg^{-1}$. In this section $K=H$, $\eta=\psi$, and $Y:=G/(H\times H)$. In particular, $\xi_1=1$.

\begin{lemma}\label{dimension-normalization}
For every object $S\in\C$, regarded as a biequivariant sheaf, we have
\[
\FPdim_\C(S)=\frac{\dim_k S}{|H|}.
\]
\end{lemma}
\begin{proof}
Let $\M:=\M(H,\psi)$. Its simple objects $S_{\bar g}$ have underlying dimension $|H|$, by Corollary \ref{simprojhpsi}. The invertible objects $\delta_g$ of $\Coh(G,\omega)$ act transitively on them. If $F$ is an exact module endofunctor of $\M$, its module constraint gives
\[
F(S_{\bar g})\cong\delta_g\ot F(S_{\bar1}).
\]
Their dimensions are therefore equal. Additivity on composition series implies
\[
\dim_kF(T)=a(F)\dim_kT\quad(T\in\M),\qquad
 a(F):=\frac{\dim_kF(S_{\bar1})}{|H|}.
\]
The function $a$ is additive, multiplicative under composition, equals $1$ on the identity, and is positive on nonzero functors. It is thus the positive dimension character on the Grothendieck ring, namely Frobenius--Perron dimension \cite[Proposition~3.3.6(3)]{EGNO}. Opposite composition does not change this scalar character. We can take $S_{\bar1}=\O(H)_\psi$, the regular comodule, which is the unique simple on the free orbit $H$. Under Lemma \ref{newsimple}, $F_S(\O(H)_\psi)=S$, hence substituting this equality gives the formula.
\end{proof}

\begin{theorem}\label{simpobjs}
The following hold.
\begin{enumerate}
\item For every $Z\in Y(k)$ with representative $g\in Z(k)$, there is an equivalence
\[
{\rm Ind}_Z:{\rm Rep}(H^g,\xi_g)_k\xrightarrow{\simeq}\C_Z.
\]
The category $\C_H$ is a tensor subcategory of $\C$, and ${\rm Ind}_H$ can be equipped with a tensor structure giving a tensor equivalence ${\rm Rep}(H)_k\simeq\C_H$.
\item The simple objects of $\C$ are parametrized by the pairs $(Z,V)$ with $V$ a simple stabilizer representation, and
\[
\C=\bigoplus_{Z\in Y(k)}\overline{\C_Z}.
\]
The Serre closure $\overline{\C_H}$ is a tensor subcategory as well.
\item For any $V\in{\rm Rep}(H^g,\xi_g)_k$,
\[
\FPdim({\rm Ind}_Z(V))=\frac{|H|}{|H^g|}\dim_kV.
\]
\item Duality sends the simple objects labelled by $Z$ to those labelled by $Z^{-1}$, and gives a contravariant equivalence
\[
\overline{\C_Z}\xrightarrow{\simeq}\overline{\C_{Z^{-1}}}.
\]
For a simple $V\in\Rep(H^g,\xi_g)_k$, use $g^{-1}$ as the representative of the inverse orbit and set
\begin{equation}\label{transported-dual-label}
V^{\vee,g}:=\operatorname{Res}_{g^{-1}}\bigl(\Ind_{Z,g}(V)^*\bigr)
\in\Rep(H^{g^{-1}},\xi_{g^{-1}})_k.
\end{equation}
Here $\operatorname{Res}_{g^{-1}}$ is fiber evaluation with its stabilizer transport, and the second subscript records the representative used in transitive descent. Then
\[
\Ind_{Z,g}(V)^*\cong\Ind_{Z^{-1},g^{-1}}(V^{\vee,g}).
\]
For another representative of $Z^{-1}$, change this label by Proposition~\ref{representative-change}.

\end{enumerate}
\end{theorem}

\begin{proof}
The equivalences and the simple-object decomposition are Theorem~\ref{simmodrep}. For the tensor statement about $\C_H$, work entirely on the subgroup $H$. Since $\omega_{\mid H}=d\psi$, gauge by $\psi$ identifies this local dual with the dual of $\Coh(H)$ relative to ${\rm Vec}$, which is ${\rm Rep}(H)$. The gauge change is trivial on the diagonal stabilizer at $1$ by normalization, so its fiber equivalence is the one used in ${\rm Ind}_H$. Extension by zero from $H$ respects the cotensor product and its unit $\O(H)_\psi$, so gives the asserted tensor subcategory. The free functor $\mathbf F_H$ (\ref{boldfzg1}) can be made tensor by transporting this tensor structure.

For (3), transitive descent gives $\dim_k{\rm Ind}_Z(V)=|Z|\dim_kV$, and $|Z|=|H|^2/|H^g|$, so the claim follows from  Lemma \ref{dimension-normalization}.

For (4), the dual of a simple is again simple. If $S$ is supported on $Z(k)$, and its dual on $W(k)$, the nonzero coevaluation $\one\to S\ot S^*$, and the support of the cotensor product imply that $1\in Z(k)W(k)$. Thus, $W=Z^{-1}$. Indeed, the cotensor product is an underlying subobject of the convolution product, whose support is contained in the product of the two supports. This proves the assertion for simples, and exactness of duality proves the assertion for Serre closures. Then \eqref{transported-dual-label} makes sense because every simple in the inverse summand is annihilated by its orbit ideal. Finally, the product of two supports lying in $H(k)$ lies in $H(k)$, so $\overline{\C_H}$ is closed under tensor product, and duality gives its rigidity.
\end{proof}

\begin{remark}\label{duality-conventions}
To express a contragredient in a chosen matrix realization, one must first transport along
\[
H^{g^{-1}}\longrightarrow H^g,\qquad h\longmapsto ghg^{-1},
\]
and match the resulting multiplier with $\xi_{g^{-1}}$ by the appropriate scalar cochain. Gauge and representative changes are those of Propositions~\ref{gauge-covariance} and \ref{representative-change}. \qed
\end{remark}

\begin{theorem}\label{projsimpobjs}
The following hold.
\begin{enumerate}
\item For a simple $V\in{\rm Rep}(H^g,\xi_g)_k$, the projective cover is the lifted-idempotent object of Proposition \ref{closedptproj}, and as an underlying sheaf
\[
P_\C(\mathbf F_Z(V))\cong
\O(X_Z)\ot_kP_{(H^g,\xi_g)}(V).
\]
Its coactions are those of \eqref{projective-coaction}, transported to this sheaf.
\item Its Frobenius--Perron dimension is
\[
\FPdim P_\C(\mathbf F_Z(V))=
\frac{|G^{\circ}||H(k)|}{|H^{\circ}||H^g(k)|}
\dim_kP_{(H^g,\xi_g)}(V).
\]
\item If the dimension of a Serre summand is defined by its regular weighted sum, then
\[
\FPdim\overline{\C_Z}:=
\sum_{S\text{ simple in }\C_Z}\FPdim(S)\FPdim(P_\C(S))
=|G^{\circ}|\,|Z(k)|.
\]
\end{enumerate}
\end{theorem}

\begin{proof}
Parts (1) and (2) follow from Theorem \ref{prsimmodrep} and Lemma \ref{dimension-normalization}. 

For (3), the algebra $A_Z$ has dimension $|X_Z||H|^2$. Over the algebraically closed field $k$, its left regular module is
\[
A_Z\cong\bigoplus_{S\text{ simple in }\C_Z}P_\C(S)^{\oplus\dim_kS}.
\]
Dividing dimensions by $|H|^2$ and using Lemma~\ref{dimension-normalization} on both factors gives the asserted weighted sum $|X_Z|=|G^{\circ}||Z(k)|$.
\end{proof}

\begin{remark}\label{unimodularity-caution}
Theorem~\ref{thm:criterion} gives the exact unimodularity criterion. In particular, unimodularity of $\Rep_k(H)$ alone is insufficient, as Example~\ref{borel-counterexample} shows. The additional character is obtained from the equivariant structure on the infinitesimal neighborhood of the identity orbit. \qed
\end{remark}

\begin{corollary}\label{semisimplicity}
Let $\C:=\C(G,\omega,H,\psi)$. The following are equivalent:
\begin{enumerate}
\item $\C$ is fusion.
\item $G^\circ\subseteq H$, and $H$ is linearly reductive.
\item $G/H$ is \'etale and $\Rep_k(H)$ is fusion.
\end{enumerate}
Thus semisimplicity, like unimodularity, is independent of the admissible twisting for a fixed inclusion $H\subseteq G$.
\end{corollary}

\begin{proof}
By Theorem \ref{projsimpobjs},
\[
\FPdim P_\C(\one)
=\frac{|G^\circ|}{|H^\circ|}\,\dim_k P_H(k).
\]
Both factors are positive integers. Since a finite tensor category is fusion if and only if its unit is projective \cite[Corollary~4.2.13]{EGNO}, semisimplicity implies that the displayed product is one, hence $G^\circ=H^\circ$, and $P_H(k)=k$.

The equality $G^\circ=H^\circ$ is equivalent to $G^\circ\subseteq H$. The condition $P_H(k)=k$ is equivalent to semisimplicity of $\Rep_k(H)$; this is linear reductivity of $H$. Thus, these two conditions give $\FPdim P_\C(\one)=1$. Since its head is $\one$, and every nonzero object has positive Frobenius--Perron dimension, $P_\C(\one)=\one$, so $\C$ is semisimple.

Finally, the local length of the homogeneous scheme $G/H$ at any closed point is $|G^\circ|/|H^\circ|$. Thus, $G/H$ is \'etale exactly when this ratio is one. Equivalently, when $G^\circ\subseteq H$, the quotient is the finite reduced homogeneous space $G(k)/H(k)$.
\end{proof}

\subsection{The \'etale case}\label{sec:The etale case}
Assume $G=G(k)$. Every orbit is open and closed, so $X_Z=Z$ and $I_Z=0$.

\begin{corollary}\label{projcoveretaleg}
The following hold.
\begin{enumerate}
\item For every simple stabilizer representation $V$,
\[
P_\C(\mathbf F_Z(V))\cong
\mathbf F_Z\bigl(P_{(H^g,\xi_g)}(V)\bigr).
\]
\item We have an equivalence of abelian categories
\[
\C\simeq\bigoplus_{Z\in Y(k)}{\rm Rep}(H^g,\xi_g)_k.
\]
In particular $\C_Z=\overline{\C_Z}$.
\item The category $\C$ is fusion if and only if ${\rm char}(k)$ does not divide $|H|$.
\end{enumerate}
\end{corollary}

\begin{proof}
The first two claims follow because the orbit ideals vanish on their summands. Part (3) follows from Corollary~\ref{semisimplicity} and Maschke's theorem for the finite group $H$. In characteristic zero the condition is automatic, recovering the fusion-case descriptions of \cite{GN,O}.
\end{proof}

\subsection{The connected case}\label{sec:The connected case}
Assume $G=G^{\circ}$. There is one closed point, and one orbit of closed points. The orbit subscheme for the $H\times H$-action through $1$ is $H$, which need not be all of $G$.

\begin{theorem}\label{projcoverconnggst}
The following hold.
\begin{enumerate}
\item $\C=\overline{\C_H}$, and $\C_H\simeq{\rm Rep}(H)_k$ as tensor categories. The bar cannot in general be omitted.
\item The simple objects are $\mathbf F_H(V)$ with $V$ simple in ${\rm Rep}(H)_k$, and their dimensions are $\FPdim\mathbf F_H(V)=\dim_kV$.
\item For each simple $V$,
\[
P_\C(\mathbf F_H(V))\cong\O(G)\ot_k P_H(V)
\]
as $\O(G)$-modules, and
\[
\FPdim P_\C(\mathbf F_H(V))=\frac{|G|}{|H|}\dim_kP_H(V).
\]
The structures are interpreted as in Remark~\ref{free-sheaf-convention}.
\end{enumerate}
\end{theorem}

\begin{proof}
These follow from Theorems \ref{simpobjs} and \ref{projsimpobjs} for $H^1=H$, $\xi_1=1$, and $X_H=G$.
\end{proof}

\begin{example}\label{restlie case}
Let $\mathfrak h\subseteq\mathfrak g$ be finite-dimensional restricted $p$-Lie algebras, and let $H\subseteq G$ be the associated height-one group schemes. Theorem \ref{projcoverconnggst} identifies the simple objects of $\C$ with those of $u(\mathfrak h)\text{-}{\rm mod}$. If $H$ is unipotent, $\C$ has the unique simple object $\one$, and
\[
P_\C(\one)\cong\O(G)\ot_kP_H(k)
\]
as an underlying sheaf. In this situation $\C$ is unimodular, since its distinguished invertible object must be its only simple object. However, this does not imply $\C\simeq\Coh(G,\omega)$ as abelian categories.

For example, take $p>2$, let $\mathfrak g$ be the restricted Heisenberg Lie algebra with $[x,y]=z$, $z$ central, and $x^{[p]}=y^{[p]}=z^{[p]}=0$, and take $G=H$, $\omega=\psi=1$. Then $\C\simeq u(\mathfrak g)\text{-}{\rm mod}$, whereas $\Coh(G)=u(\mathfrak g)^*\text{-}{\rm mod}$. Both algebras are local and basic. For the first, the central nilpotent $z$ acts by zero on every simple, and the remaining commuting nilpotent $x,y$ give only the trivial simple. The second is the local commutative coordinate algebra of $G$. The first algebra is noncommutative, and the second commutative, so their basic algebras are not isomorphic. Since Morita equivalent basic finite-dimensional algebras are isomorphic, the two abelian categories are not equivalent. \qed
\end{example}

\subsection{The normal case}\label{sec:The normal case}
Assume $H$ is normal in $G$. Then 
\[Q:=G/H=H\backslash G\]
is a finite group scheme, and $\pi:G\to Q$ is a group scheme morphism. The double cosets are the cosets $gH$, and every stabilizer $H^g$ is $H$. The cocycle $\xi_g$ need not be a coboundary.

\begin{theorem}\label{projcovernormalg}
The following hold.
\begin{enumerate}
\item For $q\in Q(k)$, and $g\in\pi^{-1}(q)(k)$,
\[
\C_{gH}\simeq{\rm Rep}(H,\xi_g)_k,
\qquad \C=\bigoplus_{q\in Q(k)}\overline{\C_{gH}}.
\]
\item If $Q$ is \'etale, these are the full summands; that is,
\[
\C\simeq\bigoplus_{q\in Q(k)}{\rm Rep}(H,\xi_g)_k.
\]
If in addition every $[\xi_g]$ is trivial, choices of trivializations give a noncanonical abelian equivalence
\[
\C\simeq\Coh(Q)\boxtimes{\rm Rep}(H)_k.
\]
\item Suppose $Q$ is \'etale, $\psi=1$, and $\omega=\pi^{*\ot3}\Omega$ for a normalized $3$-cocycle $\Omega$ on $Q$. Then there is a fully faithful tensor functor
\[
\Coh(Q,\Omega)\longrightarrow\C,
\qquad \delta_q\longmapsto\O(\pi^{-1}(q)),
\]
with the natural left and right $H$-equivariance. Its unit is $\O(H)$ in $\C$, not $\delta_1$ in $\Coh(G)$.
\end{enumerate}
\end{theorem}

\begin{proof}
Part (1) is Theorem \ref{simpobjs}. When $Q$ is \'etale, each fiber $gH$ is open and closed in $G$, hence its orbit ideal vanishes on its Serre summand. This gives the first assertion of (2), and changing a coboundary multiplier to $1$ gives the second. No tensor structure is claimed for this abelian product.

For (3), the restricted multiplier \eqref{intro-tau} is $1$, since each of its $H$-arguments maps to the identity of $Q$. The sheaf on the $H$-bitorsor $\pi^{-1}(q)$ has the stated ordinary biequivariant structure. Multiplication and descent identify the cotensor product of the objects for $q,r$ with that for $qr$. The associativity of these identifications is exactly multiplication by $\Omega(q,r,s)$, inherited from $\omega$. Their unit is the regular bicomodule $\O(H)$. Each object is simple and invertible, different cosets have no morphisms between them, and its endomorphism ring is $k$ by transitive descent. Extending on direct sums gives the fully faithful tensor functor.
\end{proof}

\begin{remark}\label{normal-obstruction}
More generally, to obtain a pointed tensor subcategory mapping onto an \'etale $Q$, one needs invertible objects $U_q\in\C_{gH}$ with $U_q\ot U_r\cong U_{qr}$. The chosen product isomorphisms then determine a $3$-cocycle $\Omega$ by their associativity defect, and give an embedding of $\Coh(Q,\Omega)$. An embedding with trivial quotient associator requires $[\Omega]=1$. Existence of each $U_q$ separately, or normality of $H$ alone, does not supply the product isomorphisms, or this last condition. Also, ordinary sheaf pullback $\pi^*:\Coh(Q)\to\Coh(G)$ is not a unital tensor functor for convolution when $H\ne1$, since it sends $\delta_1$ to $\O(H)$. \qed
\end{remark}

\begin{example}\label{normal-rank-counterexample}
Assume ${\rm char}(k)\ne2$, and take
\[
\begin{gathered}
G:=(\mathbb{Z}/2\mathbb{Z})^3,\qquad H:=(\mathbb{Z}/2\mathbb{Z})^2\times\{0\},\qquad\psi=1,\\
\omega(x,y,z)=(-1)^{x_3y_1z_2}.
\end{gathered}
\]
Trilinearity proves the cocycle identity and normalization, and $\omega_{\mid H}=1$. The subgroup $H$ is central. The identity coset has cocycle $1$, hence four one-dimensional simple objects. For $g:=(0,0,1)$, \eqref{defnrmwg} gives
\[
\xi_g(l,m)=(-1)^{l_1m_2}.
\]
The associated twisted group algebra is $M_2(k)$; its generators square to $1$ and anticommute. This coset therefore contributes one simple object of dimension $2$. Thus, $\C$ has five simple objects, with squared dimensions summing to $4+4=8$. In particular, $\C$ is not equivalent to the abelian Deligne product $\Coh(G/H)\boxtimes{\rm Rep}(H)$.

Even for $H=1$, the untwisted quotient associator cannot simply be imposed. On $G:=\mathbb{Z}/2\mathbb{Z}$, the cocycle $\omega(a,b,c)=(-1)^{abc}$ is not a coboundary, since for every normalized $2$-cochain $J$, $dJ(s,s,s)=1$, whereas $\omega(s,s,s)=-1$. Hence, the identity underlying functor from $\Coh(\mathbb{Z}/2\mathbb{Z})$ to $\Coh(\mathbb{Z}/2\mathbb{Z},\omega)$ has no tensor structure. \qed
\end{example}

\begin{example}\label{normal-untwisted-counterexample}
Assume $\operatorname{char}(k)\ne2$. Let
\[
H:=(\mathbb Z/2\mathbb Z)^4,
\qquad \sigma(a_1,a_2,a_3,a_4)=(a_3,a_4,a_1,a_2),
\]
let $G:=H\rtimes\langle s\rangle$ with $s^2=1$, and conjugation by $s$ equal to $\sigma$, and let
\[
\omega=1,\qquad \psi(a,b)=(-1)^{a_1b_2}.
\]
Bilinearity shows that $\psi$ is a normalized $2$-cocycle. There are two $H$-double cosets. On the identity coset, $\xi_1=1$, so there are sixteen one-dimensional simple objects. On the other coset, represented by $s$, \eqref{defnrmwg} gives
\[
\xi_s(a,b)=\psi(a,b)\psi(\sigma(a),\sigma(b))^{-1}
=(-1)^{a_1b_2+a_3b_4}.
\]
Its twisted group algebra is the tensor product of two copies of $M_2(k)$: in each coordinate pair the two standard generators square to one and anticommute, and generators belonging to different pairs commute. Hence, it is $M_4(k)$, and has one simple module, of dimension four.

Thus, $\C(G,1,H,\psi)$ has seventeen simple objects, with squared dimensions $16+4^2=32=|G|$. On the other hand,
$\Coh(G/H)\boxtimes\Rep(H)$ has thirty-two simple objects. Thus, the two abelian categories are not equivalent. Moreover, the nonidentity coset contains no invertible object, so a pointed tensor subcategory representing both quotient degrees cannot exist. This disproves the unrestricted normal-subgroup product assertion of \cite[Theorem~6.5(2)]{GS}, already for $\omega=1$. The extra cocycle-triviality hypothesis in Theorem~\ref{projcovernormalg}(2) cannot be omitted in general. \qed
\end{example}

\subsection{Unimodularity}\label{sec:Unimodularity}
For a finite group scheme $K$, let
$k[K]:=\O(K)^*$. Recall that its modular character
\[
\delta_K\in\Hom_{\mathrm{grp.sch.}}(K,\Gm)
\subseteq\O(K)^\times
\]
is defined by choosing a nonzero left integral $\Lambda_K\in k[K]$, and requiring
\begin{equation}\label{eq:modular-convention}
 a\Lambda_K=\eps(a)\Lambda_K,\qquad
 \Lambda_Ka=\delta_K(a)\Lambda_K\quad(a\in k[K]).
\end{equation}

We also write $\delta_K$ for the corresponding algebra homomorphism $k[K]\to k$. Products of characters correspond to convolution of these homomorphisms. Thus, for the inclusion $\iota_*:k[H]\hookrightarrow k[G]$, the equality $\delta_G|_H=\delta_H^2$ means
\[
\delta_G(\iota_*(a))
=\sum\delta_H(a_{(1)})\delta_H(a_{(2)})
\qquad(a\in k[H]).
\]
It is not enough to check this equality only on $H(k)$.
The integral exists and is unique up to a nonzero scalar; its relation with the distinguished invertible object is given in \cite[Propositions~6.5.3 and 6.5.5]{EGNO}.

\begin{lemma}\label{lem:two-socles}
Let $K$ be a finite group scheme.
\begin{enumerate}
\item The socle of the projective cover $P_{K}(k)$ of the trivial left
$k[K]$-module is
\[
 \soc P_{K}(k)\cong k_{\delta_K^{-1}}.
\]
\item Let $\O(K)_1$ be the local factor at the identity. Under the
left conjugation action on functions
\begin{equation}\label{eq:conjugation-functions}
 (h\cdot f)(g):=f(h^{-1}gh),
\end{equation}
the one-dimensional space $\soc\O(K)_1$ has character $\delta_K$.
\end{enumerate}
\end{lemma}

\begin{proof}
For (1), choose a primitive idempotent $e\in k[K]$,
whose simple head is $k_{\delta_K}$. Then $\delta_K(e)=1$, so
$\Lambda_Ke=\Lambda_K$. Thus, the trivial module $k\Lambda_K$
is a submodule of $k[K]e=P(k_{\delta_K})$.
Since a finite-dimensional Hopf algebra is Frobenius, an indecomposable
projective is also an indecomposable injective, and has simple socle.
Thus, $\soc P(k_{\delta_K})=k$. Tensoring with the invertible
module $k_{\delta_K^{-1}}$ preserves projective covers and socles, and
gives (1). 

For (2), a nonzero integral $\lambda_K\in\O(K)$ is supported at the identity
and spans $\soc\O(K)_1$, since
$f \lambda_K=\eps(f)\lambda_K$ for every $f\in\O(K)$, and the integral space is
$1$-dimensional. By the nondegenerate integral pairing for finite-dimensional Hopf algebras, we can normalize $\lambda_K$, so that
$\langle \lambda_K,\Lambda_K\rangle=1$.

The conjugation action preserves the integral line. For every
commutative $k$-algebra $T$ and $h\in K(T)$, view $h$ as a grouplike
element of $k[K]\ot_kT$. Then
$h^{-1}\Lambda_Kh=\delta_K(h)\Lambda_K$.
Pairing with $\lambda_K$ gives
\[
 \langle h\cdot \lambda_K,\Lambda_K\rangle
 =\langle \lambda_K,h^{-1}\Lambda_Kh\rangle=\delta_K(h).
\]
Therefore, $h\cdot \lambda_K=\delta_K(h)\lambda_K$. Since this holds after every
base change $k\to T$, it proves the asserted group scheme character
identity.
\end{proof}

Let $X\subseteq G$ be the open-and-closed subscheme, whose closed points
are $H(k)$; thus, $X$ contains the full connected components meeting
$H$. It is invariant under $H\times H$. Let
\begin{equation}\label{eq:relative-data}
 R:=\O(X),\,\, S:=\O(H),\,\,
 I:=\ker(R\rightarrow S),\,\, J:={\rm Ann}_R(I).
\end{equation}
The ideal $I$ is invariant and nilpotent. The module $J$ is an ordinary,
untwisted equivariant $S$-module, with equivariance induced from that
of $R$.

\begin{lemma}\label{lem:relative-line}
The $S$-module $J$ is invertible. With the left stabilizer convention,
its fiber $J(1):=k\ot_SJ$ at the identity satisfies
\begin{equation}\label{eq:relative-character}
 J(1)\cong k_{\jmath_{G,H}},\qquad
 \jmath_{G,H}:=(\delta_G|_H)\delta_H^{-1}.
\end{equation}
\end{lemma}

\begin{proof}
The commutative algebras $R$ and $S$ are Frobenius, since $S$ is a finite-dimensional 
Hopf algebra, and $R$ is a direct factor of the finite-dimensional Hopf algebra
$\O(G)$. Therefore, as $S$-modules,
\[
 J=\Hom_R(S,R)\cong\Hom_R(S,R^*)\cong S^*\cong S.
\]
Thus, $J$ is invertible.

Let $R_1:=\O(G)_1$, and $S_1:=\O(H)_1$ be the local factors at the
identity, and let $J_1$ be the corresponding factor of $J$.
Multiplication gives a canonical equivariant isomorphism
\begin{equation}\label{eq:relative-socle}
 \soc(S_1)\ot_k J(1)\xrightarrow{\simeq}\soc(R_1).
\end{equation}
To see this, $J_1$ is free of rank one over $S_1$, and
\[
 \soc_{S_1}(J_1)=\soc_{R_1}(R_1).
\]
Indeed, the inverse image of the maximal ideal of $S_1$ is the maximal
ideal of $R_1$, and every element of $\soc(R_1)$ is annihilated by $I$.
For a free rank-one module, multiplication identifies its socle with
the socle of the ring tensored with its residue-field fiber. This
proves \eqref{eq:relative-socle}, including its equivariance.

On the diagonal stabilizer the right action on the base is $x\mapsto h^{-1}xh$. Ordinary forward transport induces inverse pullback on sections. Converting this to a left representation therefore gives the action \eqref{eq:conjugation-functions}. By Lemma~\ref{lem:two-socles}(2), the characters of $\soc(R_1)$ and $\soc(S_1)$ are $\delta_G|_H$ and $\delta_H$, respectively. Their quotient in \eqref{eq:relative-socle} is \eqref{eq:relative-character}.

\end{proof}

Let
\[
 \E_H:\C_H\xrightarrow{\simeq}\Rep_k(H)
\]
be restriction to the fiber at $1$, with the stabilizer action converted
to a \emph{left} action by inversion. Thus, $h$ acts on this fiber by the
transport for $(h^{-1},h^{-1})$. Write ${\bf F}_H^{\ell}$ for its quasi-inverse.
The tensor unit has trivial stabilizer action, so
$\E_H(\one)=k$, and ${\bf F}_H^{\ell}(k)=\one$.

\begin{theorem}\label{thm:criterion}
The following conditions are equivalent:
\begin{enumerate}
\item $\C(G,\omega,H,\psi)$ is unimodular.
\item $\soc P_{\C}(\one)\cong\one$.
\item $\delta_G|_H=\delta_H^2$ as characters of the group scheme $H$.
\end{enumerate} 
Thus, for a fixed inclusion $H\subseteq G$, unimodularity is independent
of $\omega$ and $\psi$.
\end{theorem}

\begin{proof}
The equivalence of (1) and (2) is the standard characterization of
unimodularity by the socle of the projective cover of the unit;
see \cite[Sections~6.4--6.5]{EGNO}. 

We prove
\begin{equation}\label{eq:socle-formula}
 \soc P_{\C}(\one)
 \cong{\bf F}_H^{\ell}\bigl(k_{\chi_{G,H}}\bigr),\qquad
 \chi_{G,H}:=(\delta_G|_H)\delta_H^{-2},
\end{equation}
where $k_\chi$ denotes the one-dimensional left representation with
character $\chi$, which also proves the equivalence with (3).

Let $F:=H\times H$, and let $\Phi$ be the action-groupoid cocycle
coming from $\omega,\psi$ in the biequivariant realization. Recall that 
$A:=A(X,F,\Phi)$, $\bar A:=A/IA=A(H,F,\Phi|_H)$,
$A\text{-}\mathrm{mod}\simeq\overline{\C_H}$, and $\bar A\text{-}\mathrm{mod}\simeq\C_H$.

Let $P:=P_\C(\one)$. Since the unit belongs to the Serre summand
$\overline{\C_H}$, so does $P$. As a projective $A$-module, $P$ is a
direct summand of a finite sum of copies of $A$, and is therefore
projective over $R=\O(X)$.

Because $IA$ is nilpotent, reduction of $P$ gives the projective
cover of the unit over $\bar A$:
\begin{equation}\label{eq:reduced-projective}
 \overline{P}:=P/IP\cong P_{\C_H}(\one),\qquad
 \E_H(\overline{P})\cong P_{\Rep_k(H)}(k).
\end{equation}
One can verify this directly using a primitive idempotent:
if $P=Ae$, then $\overline{P}=\overline{A}\bar e$, and $\bar e$ remains primitive
because $IA$ is nilpotent. The head remains the unit.

Multiplication gives an isomorphism
\begin{equation}\label{eq:annihilator-projective}
 J\ot_S\bar P\xrightarrow{\simeq}{\rm Ann}_P(I),
 \qquad j\ot\bar v\longmapsto jv.
\end{equation}
It is well defined because $JI=0$. For $P=R$, it is the identity
$J\ot_SS\cong J$, and hence it is an isomorphism for finite free
$R$-modules, and their direct summands, which proves it for $P$.

Crucially, \eqref{eq:annihilator-projective} respects the twisted
equivariant structures. Transport on $P$ is semilinear over transport
on $R$, so transporting $jv$ gives the product of the transported
$j$ and $v$. The line $J$ has ordinary equivariance, while $\overline{P}$
and ${\rm Ann}_P(I)$ have the same multiplier $\Phi|_H$. Thus, multiplication
has exactly the required multiplier, with no extra scalar factor.

Every simple $A$-module is annihilated by the nilpotent ideal $IA$. Hence the simple subobjects of $P$ are precisely those of ${\rm Ann}_P(I)$, viewed as an $\overline A$-module. Thus,
\begin{equation}\label{eq:socle-ann}
 \soc_\C(P)=\soc_{\C_H}\bigl({\rm Ann}_P(I)\bigr).
\end{equation}

On the identity orbit, restriction of $\Phi$ to the diagonal
stabilizer is $1$, since the two factors involving $\psi$ cancel,
and normalization of $\omega$ makes all remaining factors equal to
$1$. Transitive descent therefore identifies $\C_H$ with ordinary
representations of $H$. Moreover, the unit $\O(H)_\psi$ has transport
\[
 \theta_{(h,h)}(1)=\omega(h,h^{-1},h)
 \frac{\psi(h,h^{-1})}{\psi(h^{-1},h)}\,\id=\id,
\]
by $d\psi=\omega|_H$. Thus, its fiber is the trivial representation.
Under $\E_H$, tensoring an equivariant sheaf
by the ordinary equivariant line $J$ corresponds to tensoring its
stabilizer representation by $J(1)$.

Applying $\E_H$ to \eqref{eq:annihilator-projective}, and using
\eqref{eq:reduced-projective}, gives
\[
 \E_H\bigl({\rm Ann}_P(I)\bigr)
 \cong J(1)\ot_kP_{\Rep_k(H)}(k).
\]
Since tensoring by a one-dimensional representation is an exact
autoequivalence, Lemmas~\ref{lem:two-socles} and
\ref{lem:relative-line}, together with \eqref{eq:socle-ann}, yield
\begin{equation*}
 \E_H\bigl(\soc_\C(P)\bigr)
 \cong k_{\jmath_{G,H}}\ot_k k_{\delta_H^{-1}}\cong k_{(\delta_G|_H)\delta_H^{-2}},
\end{equation*}
which is \eqref{eq:socle-formula}. Since $\E_H(\one)=k$, the socle
is the unit if and only if $\delta_G|_H=\delta_H^2$.
\end{proof}

\begin{corollary}\label{cor:lie}
Suppose $p>0$, and let $\mathfrak h\subseteq\mathfrak g$ be finite-dimensional
restricted Lie algebras, with associated height-one group schemes
$H\subseteq G$. The following are equivalent:
\begin{enumerate}
\item
 $\C(G,\omega,H,\psi)$ is unimodular.
 \item 
 $\tr(\operatorname{ad}_{\mathfrak g}x)
 =2\tr(\operatorname{ad}_{\mathfrak h}x)$ 
for every $x\in\mathfrak h$.
\item
$\tr\bigl(x\mid\mathfrak g/\mathfrak h\bigr)
 =\tr(\operatorname{ad}_{\mathfrak h}x)$, 
$x\in\mathfrak h$,
where the action on $\mathfrak g/\mathfrak h$ is induced by the adjoint
action.
\end{enumerate}
\end{corollary}

\begin{proof}
We have $k[G]=u(\mathfrak g)$, and $k[H]=u(\mathfrak h)$. With convention
\eqref{eq:modular-convention},
\begin{equation}\label{eq:lie-modular}
 \delta_G(x)=\tr(\operatorname{ad}_{\mathfrak g}x)
 \qquad(x\in\mathfrak g).
\end{equation}
To verify the sign, choose a basis $x_1,\ldots,x_n$
of $\mathfrak g$, and use the PBW filtration. The leading symbol of
a nonzero left integral is annihilated by every $x_i$ in
\[
 \gr u(\mathfrak g)
 =k[x_1,\ldots,x_n]/(x_1^p,\ldots,x_n^p).
\]
It is therefore a nonzero multiple of
$x_1^{p-1}\cdots x_n^{p-1}$. The adjoint action of $x$ on this
one-dimensional top degree has eigenvalue
\[
 (p-1)\tr(\operatorname{ad}_{\mathfrak g}x)
 =-\tr(\operatorname{ad}_{\mathfrak g}x).
\]
On the other hand,
$[x,\Lambda_G]=-\delta_G(x)\Lambda_G$ by
\eqref{eq:modular-convention}. Thus, comparing leading symbols proves
\eqref{eq:lie-modular}. The same argument applies to $\mathfrak h$.

Since $x\in\mathfrak h$ is primitive,
$(\delta_H*\delta_H)(x)=2\delta_H(x)$. Since $u(\mathfrak h)$ is generated by $\mathfrak h$, and the characters are algebra homomorphisms, the equality $\delta_G|_H=\delta_H^2$ is equivalent to
$\operatorname{tr}(\operatorname{ad}_{\mathfrak g}x)
=2\operatorname{tr}(\operatorname{ad}_{\mathfrak h}x)$ 
$(x\in\mathfrak h)$.
This proves the equivalence of (1) and (2). 

Finally, since trace is additive on the
$\mathfrak h$-module exact sequence
\[0\to\mathfrak h\to\mathfrak g\to\mathfrak g/\mathfrak h\to0,\]
(3) is equivalent to (2).
\end{proof}

\begin{example}\label{borel-counterexample}\label{616}
Let $p>2$, and let $\mathfrak g$ have basis $t,x$ with
\[
 [t,x]=x,\qquad t^{[p]}=t,\qquad x^{[p]}=0,
 \qquad\mathfrak h:=kt.
\]
Then $G:=\alpha_p\rtimes\mu_p$, and $H:=\mu_p$. We have
$\tr(\operatorname{ad}_{\mathfrak g}t)=1$,
and $\tr(\operatorname{ad}_{\mathfrak h}t)=0$.
Thus, by Corollary~\ref{cor:lie}, $\C(G,1,H,1)$ is not
unimodular, although $\Rep_k(H)$ is semisimple and hence unimodular. Thus the implication in \cite[Theorem~6.2(5)]{GS} is false in general; Theorem~\ref{thm:criterion} supplies the missing modular condition. \qed
\end{example}

\begin{example}\label{modular-cancellation}
Let $p>2$, let $\mathfrak g$ be the restricted semidirect product with
basis $t,x,y$, and
\[
 [t,x]=x,\qquad [t,y]=y,\qquad[x,y]=0,
 \qquad t^{[p]}=t,\quad x^{[p]}=y^{[p]}=0,
\]
and let $\mathfrak h:=kt\oplus kx$. Then
$\tr(\operatorname{ad}_{\mathfrak g}t)=2$, $\tr(\operatorname{ad}_{\mathfrak h}t)=1$,
and the two traces at $x$ are zero. Thus,
$\C(G,\omega,H,\psi)$ is unimodular for every admissible pair, whereas
$\Rep_k(H)$ is not unimodular. The two modular contributions need
not vanish separately; they can cancel. \qed
\end{example}

\begin{corollary}\label{cor:checks}
The following hold.
\begin{enumerate}
\item If $G$ is \'etale, then $\C(G,\omega,H,\psi)$ is unimodular,
even when it is not semisimple.
\item If $H$ has no nontrivial characters, in particular if $H$ is
unipotent, then $\C(G,\omega,H,\psi)$ is unimodular.
\item If $\Rep_k(G)$ and $\Rep_k(H)$ are unimodular, then so is $\C$.
\item If $\Rep_k(G)$ is unimodular, then
$\C(G,\omega,H,\psi)$ is unimodular if and only if $\delta_H^2=1$.
\item For $H=1$, the criterion is automatic. For $H=G$, it is
$\delta_G=1$.
\end{enumerate}
\end{corollary}

\begin{proof}
If $G$ is \'etale, both $k[G]$ and $k[H]$ are group algebras. The sum
of the group elements is a two-sided integral, so
$\delta_G=\delta_H=1$, in every characteristic. If $H$ has no
nontrivial characters, both characters in Theorem \ref{thm:criterion}(3) are
trivial. A finite unipotent group scheme has only the trivial simple
representation, and hence has no nontrivial characters. The remaining
claims are immediate from Theorem~\ref{thm:criterion}.
\end{proof}

\begin{example}
Let $K$ be a finite group scheme, $G:=K\times K$, and
$H:=\Delta(K)$. A tensor product of left integrals shows that
\[
 \delta_{K\times K}=\delta_K\boxtimes\delta_K,
 \qquad\delta_{K\times K}|_{\Delta(K)}=\delta_K^2.
\]
Thus, the criterion is automatic for the diagonal subgroup, for every
admissible twisting. In particular, in the untwisted case it recovers
unimodularity of the center of $\Coh(K)$. \qed
\end{example}

\section{\texorpdfstring{The center of $\Coh(G,\omega)$}{The center of Coh(G, omega)}}\label{sec:twisted-double}

Fix $G,\omega$ as in \S\ref{sec:Coh(G,w)}, and let $\D:=\Coh(G,\omega)$. An object of its center $\Z(G,\omega):=\Z(\D)$ is a pair $(M,c)$, where $M\in\D$ and
\[
c_Y:Y\ot M\xrightarrow{\cong}M\ot Y\qquad(Y\in\D)
\]
is a natural half-braiding. The hexagon uses the associator of \S\ref{sec:Conventions}. The center is a finite braided tensor category, of Frobenius--Perron dimension $|G|^2$. Its braiding is nondegenerate in the finite tensor category sense; see \cite{EGNO,Sh}.

\subsection{The conjugation description}\label{sec:conjugation-description}
By \cite[Proposition~7.13.8]{EGNO}, the center is the dual of $\D\boxtimes\D^{\mathrm{rev}}$ with respect to the regular bimodule category $\D$, up to the dual convention of \S\ref{sec:Conventions}. For the orbit and projective descriptions we use instead a direct realization by twisted conjugation-equivariant sheaves. 

Consider the right conjugation action $x\cdot a=a^{-1}xa$, and define 
\begin{equation}\label{conjugation-multiplier}
\Phi_{\rm ad}(x;a,b):=
\frac{\omega(x,a,b)\omega(a,b,(ab)^{-1}x(ab))}
{\omega(a,a^{-1}xa,b)}.
\end{equation}
The same four pentagon substitutions as in Proposition \ref{taulemma}, now with $h_i=k_i$, show that this is a normalized action-groupoid cocycle (the two remaining factors in \eqref{omega-factor-ratio} cancel).

\begin{proposition}\label{simpobjscenter}
There is an equivalence
\[
\Z(G,\omega)\simeq\Coh^{(G,\Phi_{\rm ad})}(G),
\]
preserving underlying sheaves. The tensor product on the right is the transported convolution product, with the associators in the half-braiding formula. We also denote this category by $\Coh^{(G)}(G,\omega)$.
\end{proposition}

\begin{proof}
Naturality on the regular $\O(G)$-module shows that a half-braiding is determined by
\[
\rho(m)=c_{\O(G)}(1\ot m)=\sum m_0\ot m_1.
\]
Indeed, for any $Y\in\Coh(G)$, and $y\in Y$, the module map $\O(G)\to Y$, $f\mapsto fy$, gives
$c_Y(y\ot m)=\sum m_0\ot m_1y$.
Linearity over $\O(G)$ is exactly compatibility with conjugation, and the unit constraint is counitality.

Write $C_a(x):M_x\to M_{axa^{-1}}$ for the fiber map induced by $c$. The half-braiding hexagon is
\begin{equation}\label{half-braiding-hexagon}
c_{A\ot B}
=\alpha_{M,A,B}(c_A\ot\id_B)\alpha^{-1}_{A,M,B}
(\id_A\ot c_B)\alpha_{A,B,M}.
\end{equation}
On the pulled-back fibers it gives
\[
C_{ab}(x)=K(x;a,b)\,C_a(bxb^{-1})C_b(x),
\]
where
\[
K(x;a,b):=
\frac{\omega(a,b,x)\omega(abx(ab)^{-1},a,b)}
{\omega(a,bxb^{-1},b)}.
\]
Define forward transport for right conjugation by
\[
\theta_a(y):=C_a(a^{-1}ya)^{-1}:M_y\longrightarrow M_{a^{-1}ya}.
\]
Applying the preceding hexagon identity at $x=(ab)^{-1}y(ab)$, and taking inverses, gives 
\[
\theta_b(a^{-1}ya)\theta_a(y)
=K((ab)^{-1}y(ab);a,b)\,\theta_{ab}(y).
\]
The scalar on the right is exactly $\Phi_{\mathrm{ad}}(y;a,b)$ in \eqref{conjugation-multiplier}. The calculation is on the scheme of composable arrows, so it also applies to nonreduced group schemes.

Conversely, a transport satisfying this identity gives the map $\rho$, and hence all the displayed $c_Y$. Its inverse-arrow identity ensures invertibility, and the groupoid identity gives the half-braiding hexagon. 

Finally, these constructions are inverse on objects and morphisms. The usual tensor product of half-braidings, with the same associators, gives the asserted transported tensor structure.
\end{proof}

\subsection{\texorpdfstring{The structure of $\Coh^{(G)}(G,\omega)$}{The structure of conjugation-equivariant sheaves}}
Let ${\rm C}$ be the finite quotient for conjugation. For $C\in{\rm C}(k)$, choose $g\in C(k)$, and write $G_C:=C_G(g)$ for the scheme-theoretic centralizer of $g$. Let $X_C\subseteq G$ be the union of the full local components at $C(k)$. The category $\Z_C$ consists of objects annihilated by the orbit ideal, whereas $\overline{\Z_C}$ consists of all objects supported on $X_C$.

\begin{lemma}\label{omegagprop}
The restriction of \eqref{conjugation-multiplier} to $\{g\}\times G_C^2$ is $\omega_g|_{G_C^2}$, with $\omega_g$ as in \eqref{omegag}. It is a normalized $2$-cocycle on $G_C$.
\end{lemma}
\begin{proof}
For $a,b\in G_C$, conjugation fixes $g$, so the restriction is
\[
\frac{\omega(g,a,b)\omega(a,b,g)}{\omega(a,g,b)}.
\]
This is \eqref{omegag} on $G_C$, and restriction of the groupoid identity proves the cocycle identity.
\end{proof}

\begin{theorem}\label{simpobjsnew}
Set $\Z:=\Z(G,\omega)$. The following hold.
\begin{enumerate}
\item There are equivalences
\[
\mathbf F_C:{\rm Rep}(G_C,\omega_g)_k\xrightarrow{\simeq}\Z_C,
\qquad V\longmapsto\O(C)\ot_kV,
\]
where the structure is the transitive descent of \S\ref{S:Biequivariant sheaves on group schemes}, transported to a free sheaf. These label all the simple objects, without repetitions after fixing representatives.
\item We have a decomposition into Serre summands
\[
\Z=\bigoplus_{C\in{\rm C}(k)}\overline{\Z_C}.
\]
The identity orbit category is the canonical symmetric tensor subcategory ${\rm Rep}(G)_k\simeq\Z_{\{1\}}$. Its Serre closure is generally larger.
\item For every $V$ as in (1),
\[
\FPdim\mathbf F_C(V)=|C|\dim_kV=\frac{|G|}{|G_C|}\dim_kV.
\]
For a simple $V$, its dual belongs to the inverse conjugacy orbit. Its stabilizer label is obtained by evaluating this categorical dual at $g^{-1}$, as in \eqref{transported-dual-label}, with the transported multiplier.
\item For a simple $V$, there is an isomorphism of underlying sheaves
\[
P_\Z(\mathbf F_C(V))\cong\O(X_C)\ot_kP_{(G_C,\omega_g)}(V).
\]
Its full structure is the lifted-idempotent construction for
$A(X_C,G,\Phi_{\rm ad})$, including \eqref{projective-coaction}. In particular,
\[
\FPdim P_\Z(\mathbf F_C(V))=
\frac{|G|}{|G_C(k)|}\dim_kP_{(G_C,\omega_g)}(V).
\]
\item With the regular weighted-sum definition,
\[
\FPdim\overline{\Z_C}=\frac{|G|^2}{|G_C(k)|}.
\]
\end{enumerate}
\end{theorem}

\begin{proof}
Apply transitive descent to $G_C\backslash G\cong C$, with multiplier \eqref{conjugation-multiplier}. Apply the invariant-idempotent support argument of Theorem~\ref{simmodrep} to conjugation. The identity orbit is the reduced point $1$, with stabilizer $G$ and cocycle $1$; its tensor structure and symmetric embedding follow either from the half-braiding formula, or from the usual representation subcategory of the double.

The forgetful tensor functor $\Z(\D)\to\D$ preserves Frobenius--Perron dimension, and on $\D$ this is the underlying vector-space dimension. This proves (3). 

Duality and inversion of supports give its last assertion. The proof of Theorem \ref{prsimmodrep} applies verbatim to any finite action groupoid, so it gives (4), with
\[
|X_C|=|G^{\circ}|\,|C(k)|=\frac{|G|}{|G_C(k)|}.
\]

Finally $A(X_C,G,\Phi_{\rm ad})$ has dimension $|X_C||G|$. Decomposing its regular module into projective covers, with multiplicities the dimensions of the simple modules, gives (5). 
\end{proof}

\begin{remark}\label{centerfib}
Corollary~\ref{fibfungth} applies to a group scheme-theoretical presentation of the center once its parameters are fixed. The twisted centers also occur in the classification of minimal nondegenerate extensions in \cite{G2}. For tensor subcategories and centralizers of the untwisted double, see \cite{AG}. \qed
\end{remark}

\begin{example}\label{restlie case double}
If $G$ is connected, then
$\Z(G,\omega)=\overline{{\rm Rep}(G)_k}$.
For each simple $V\in{\rm Rep}(G)_k$, Theorem \ref{simpobjsnew} gives the underlying sheaf $\O(G)\ot_kP_G(V)$ of its projective cover. It does not give the coaction $\id\ot\rho_{P_G(V)}$. The following computation shows explicitly why this distinction matters. \qed
\end{example}

\begin{lemma}\label{alpha2cocycle}
Assume ${\rm char}(k)=2$. For $\O(\alpha_2)=k[x]/(x^2)$, with $x$ primitive,
$\omega_\lambda:=1+\lambda x\ot x\ot x$ $(\lambda\in k)$ 
is a normalized $3$-cocycle. If $\lambda\ne0$, it is not a coboundary.
\end{lemma}

\begin{proof}
Use independent square-zero variables $x_1,x_2,x_3,x_4$. Each factor of the pentagon is $1$ plus a cubic monomial, or a sum of cubic monomials. Products of two such cubic terms vanish, and the identity reduces to
\[
x_1x_2(x_3+x_4)+(x_1+x_2)x_3x_4
=x_2x_3x_4+x_1(x_2+x_3)x_4+x_1x_2x_3.
\]
This equality is immediate by expansion. The inverse of $\omega_\lambda$ is itself, and normalization is clear.

Every normalized $2$-cochain is $J=1+c x_1x_2$, since its restrictions to either coordinate axis are $1$. In three square-zero variables,
\[
J(x,y+z)J(y,z)=1+c(xy+xz+yz)
=J(x+y,z)J(x,y).
\]
The possible quadratic products vanish. Thus, $dJ=1$ for every such $J$, whereas $\omega_\lambda\ne1$ for $\lambda\ne0$.
\end{proof}

\begin{theorem}\label{alpha2center}
In characteristic $2$, let $\D_\lambda:=\Coh(\alpha_2,\omega_\lambda)$. There is an equivalence of abelian categories
\begin{equation}\label{alpha2-algebra}
\Z(\D_\lambda)\simeq B_\lambda\text{-}{\rm mod},
\qquad B_\lambda:=k[X,T]/(X^2,T^2-\lambda X).
\end{equation}
Under this equivalence, $X$ is the action of $x$ on the underlying sheaf, and the half-braiding is
\begin{equation}\label{alpha2-halfbraiding}
c_Y(y\ot m)=m\ot y+Tm\ot xy.
\end{equation}
For $\lambda=0$, the algebra is $k[X,T]/(X^2,T^2)$, and for $\lambda\ne0$, it is isomorphic to $k[T]/(T^4)$.
\end{theorem}

\begin{proof}
Let $(M,c)\in\Z(\D_\lambda)$, and let $X$ be the action of $x$ on $M$, so $X^2=0$. Naturality on the regular module $A=k[x]/(x^2)$ implies that
\[
c_A(1\ot m)=m\ot1+Tm\ot x.
\]
The coefficient of $1$ is $m$ by naturality for $A\to k$, and the unit condition. Naturality for $A\to Y$, $a\mapsto ay$, gives \eqref{alpha2-halfbraiding} for all $Y$. 

Conversely, every linear $T$ gives such a natural family. The family is $A$-linear exactly when $XT=TX$, by comparing the actions of $x$ on its two sides. It is invertible; after the vector-space flip it is $1+T\ot x$, whose inverse is itself.

It remains to check the hexagon. On $U\ot V\ot M$, let $a,b,X$ denote the commuting square-zero actions of $x$ on the three factors. The three associators in \eqref{half-braiding-hexagon} each contribute $1+\lambda abX$, since the characteristic is two, and $\omega_\lambda$ is invariant under permutation of its three arguments. Their product is again $1+\lambda abX$. The two half-braidings contribute $(1+aT)(1+bT)$. Hence, the right side of the hexagon, after the common flip, is
\[
(1+\lambda abX)(1+aT)(1+bT)
=1+(a+b)T+ab(T^2+\lambda X).
\]
The left side, $c_{U\ot V}$, is $1+(a+b)T$. Thus, the hexagon holds for all $U,V$ precisely when $T^2=\lambda X$. Necessity can be checked with both $U,V$ regular, since $ab$ is then nonzero. These are exactly the relations of $B_\lambda$, and morphisms are exactly the maps commuting with $X,T$. This proves the equivalence, including its converse. For $\lambda\ne0$, eliminate $X=\lambda^{-1}T^2$; the remaining relation is $T^4=0$.
\end{proof}

\begin{corollary}\label{alpha2extensions}
The center in Theorem \ref{alpha2center} has one simple object $k$, and one indecomposable projective cover $P:=B_\lambda$. Writing $\Jac(B_\lambda)$ for the Jacobson radical and $\LL(P)$ for the Loewy length, its radical layers (from head to socle) and self-extension space are
\[
\begin{array}{c|c|c|c}
 &\text{dimensions of radical layers}&\LL(P)&\dim_k\Ext^1(k,k)\\ \hline
\lambda=0&1,\ 2,\ 1&3&2\\
\lambda\ne0&1,\ 1,\ 1,\ 1&4&1.
\end{array}
\]
In particular, the two cases are not equivalent as abelian categories.
\end{corollary}

\begin{proof}
For $\lambda=0$, the radical is $(X,T)$, its square is $kXT$, and its cube is zero. For $\lambda\ne0$, the radical is $(T)$ in $k[T]/(T^4)$. These give the radical layers. For either local augmented algebra, applying $\Hom(-,k)$ to the beginning of a projective resolution gives
\[
\Ext^1(k,k)\cong\Hom_k(\Jac(B_\lambda)/\Jac(B_\lambda)^2,k).
\]
This gives the dimensions in the table. Since an abelian equivalence preserves projective covers and extension groups, the two cases cannot be equivalent.
\end{proof}

As a module over $k[X]/(X^2)$, $P:=B_\lambda$ is free with basis $1,T$. In this basis the half-braiding operator on its projective cover is
\begin{equation}\label{alpha2-projective-matrix}
T_P:=\begin{pmatrix}0&\lambda X\\1&0\end{pmatrix},
\qquad T_P^2=\lambda X\,\id.
\end{equation}
Note that for $\lambda=1$, a putative coaction acting only on the stabilizer-projective factor would instead have $T_P^2=0$, hence would not satisfy the hexagon. Thus, \eqref{alpha2-projective-matrix} is the actual structure correction, whereas a dimension count alone could not detect the missing term.

\subsection{Restriction and relative centers}
Recall the split sequence $1\to G^{\circ}\xrightarrow{i}G\mathrel{\mathop{\rightleftarrows}^{\pi}_{q}}G(k)\to1$, and let
\begin{equation}\label{defnwcircwk}
\omega^{\circ}:=i^{*\ot3}\omega,\qquad
\omega(k):=q^{*\ot3}\omega.
\end{equation}
In general, $\pi^{*\ot3}\omega(k)\ne\omega$. Note that the splitting says $q^*\pi^*=\id$ on $\O(G(k))$, 
so one cannot recover the associator from its restriction to the reduced subgroup.
Also, the splitting maps do not automatically yield tensor functors between their centers, as the following two examples demonstrate.

\begin{example}\label{restriction-not-tensor}
Take $G:=\mathbb{Z}/2\mathbb{Z}=\{1,s\}$, $\omega=1$. In $\Z(\Coh(\mathbb{Z}/2\mathbb{Z}))$ let $M$ be the central one-dimensional object supported at $s$, with trivial character. Then
$M\ot M\cong\one$, $i^*M=0$, and $i^*(M\ot M)=k$.
Therefore $i^*$ cannot have tensor-product isomorphisms. This contradicts \cite[Theorem~8.5(4)]{GS}. The same example, with $q=\id$, shows that $i^*q_*$ takes the identity-degree part, and is not the usual forgetful functor. \qed
\end{example}

\begin{example}\label{reduced-inclusion-fails}
Let $p$ be odd, and $G:=\alpha_p\rtimes \mathbb{Z}/2\mathbb{Z}$, with $s$ acting by $t\mapsto-t$, and take $\omega=1$. On points over a test algebra with $t^p=0$, we have
$(t,1)^{-1}(0,s)(t,1)=(-2t,s)$.
Thus, the reduced subscheme $\{s\}$ is not stable under conjugation by $G$. The pushforward along $q$ of the degree-$s$ central object of $\Z(\Coh(\mathbb{Z}/2\mathbb{Z}))$ has underlying sheaf $\delta_s$, whose annihilator is not conjugation invariant. It cannot be a conjugation-equivariant sheaf on $G$. Thus the lift of $q_*$ asserted in \cite[Theorem~8.5(2)]{GS} fails even to define an object, independently of any associator. \qed
\end{example}

The valid general forgetful construction keeps the underlying sheaf on $G$, and restricts only its half-braiding. Set $\D^{\circ}:=\Coh(G^{\circ},\omega^{\circ})$. The relative center $\Z_{\D^{\circ}}(\D)$ consists of pairs $(M,c)$ with $M\in\D$, and a half-braiding $c_Y$ only for $Y\in\D^{\circ}$. There is a tensor functor
\begin{equation}\label{obvious forgetful tensor functor}
\Z(\D)\longrightarrow\Z_{\D^{\circ}}(\D),
\qquad(M,c)\longmapsto(M,c|_{\D^{\circ}}).
\end{equation}
The regular-module argument of Proposition~\ref{simpobjscenter} identifies the target with
\[
\Coh^{(G^{\circ},\Phi_{\rm ad}|_{G\times(G^{\circ})^2})}(G).
\]
It is finite abelian by Proposition~\ref{descent-algebra}. Its forgetful functor to $\D$ is faithful and exact; tensor products and duals come from relative half-braidings. Its unit is simple and its tensor product is exact in each variable, so it is a finite tensor category. The functor \eqref{obvious forgetful tensor functor} is exact and dominant by Theorem~\ref{gnntheorem} below.

\subsection{\texorpdfstring{$\Z(G,\omega)$ as $G(k)$-equivariantization}{Z(G, omega) as G(k)-equivariantization}}
Write $\Gamma:=G(k)$, and $u_a:=\delta_a$ for $a\in\Gamma$. The component decomposition gives a faithful grading
\begin{equation}\label{gamma}
\D=\bigoplus_{a\in\Gamma}\D_a,\qquad
\D_a:=\D^{\circ}\ot u_a.
\end{equation}
Tensoring by $u_a$ is an equivalence $\D^{\circ}\xrightarrow{\simeq}\D_a$. This grading retains the full associator $\omega$, including its mixed terms.

We give the relative-center argument in this setting, so that no semisimplicity or division by $|\Gamma|$ is needed. In the fusion setting this is \cite[Theorem~3.5]{GNN}.

\begin{theorem}\label{gnntheorem}
There is a coherent action of $\Gamma$ by tensor autoequivalences on $\Z_{\D^{\circ}}(\D)$, and an equivalence of tensor categories
\[
F:\Z(\D)\xrightarrow{\simeq}\Z_{\D^{\circ}}(\D)^\Gamma.
\]
The action carries $\Z_{\D^{\circ}}(\D_a)$ to $\Z_{\D^{\circ}}(\D_{hah^{-1}})$ under $h$. The functor \eqref{obvious forgetful tensor functor} identifies with forgetting equivariance; it is exact and dominant.
\end{theorem}

\begin{proof}
Choose a left dual $v_a$ of $u_a$, with evaluation $\operatorname{ev}_a:v_a\ot u_a\to\one$ and coevaluation $\operatorname{coev}_a:\one\to u_a\ot v_a$. Since $u_a$ is invertible, both maps are isomorphisms. Take $u_1=v_1=\one$ and choose normalized isomorphisms $m_{a,b}:u_a\ot u_b\to u_{ab}$. Write
\[
T_a(M):=(u_a\ot M)\ot v_a.
\]
Its tensorator $T_a(M)\ot T_a(N)\to T_a(M\ot N)$ contracts the middle $v_a\ot u_a$ by $\operatorname{ev}_a$, and its unit map is $\operatorname{coev}_a$, with the unit constraints understood. Rebracketing uses the ambient associator. The triangle identities give the unit diagrams, and the associator pentagon gives the tensor-functor diagram. In particular, $T_a$ is a tensor autoequivalence of $\D$ preserving $\D^\circ$.

Let
\[
\widehat m_{a,b}:v_b\ot v_a\longrightarrow v_{ab}
\]
be the inverse of the dual of $m_{a,b}$, under the canonical dual-of-a-tensor-product identification. Define
\begin{equation}\label{action-coherence-map}
\mu_{a,b}(M):T_aT_b(M)\longrightarrow T_{ab}(M)
\end{equation}
by $m_{a,b}$ on the left of $M$ and $\widehat m_{a,b}$ on the right, with the intervening associators. The action identity is
\begin{equation}\label{action-pentagon}
\mu_{ab,c}(M)\mu_{a,b}(T_c(M))
=\mu_{a,bc}(M)T_a(\mu_{b,c}(M)).
\end{equation}
Indeed, the two composites from $(u_a\ot u_b)\ot u_c$ to $u_{abc}$ differ by a scalar. The corresponding inverse dual composites differ by its inverse, so these discrepancies cancel on the two sides of $M$. All remaining rebracketings agree by the pentagon. Naturality of evaluation makes $\mu_{a,b}$ compatible with the tensorators, and normalization gives the action unit identities.

For $(M,c)\in\Z_{\D^\circ}(\D)$ and $Y\in\D^\circ$, transport the half-braiding by
\begin{equation}\label{relative-action-halfbraiding}
\begin{split}
Y\ot T_a(M)&\cong T_a(T_a^{-1}(Y)\ot M)
 \xrightarrow{\ T_a(c_{T_a^{-1}(Y)})\ }
 T_a(M\ot T_a^{-1}(Y))\\
&\cong T_a(M)\ot Y.
\end{split}
\end{equation}
Here $T_a^{-1}$ is a chosen monoidal quasi-inverse, and the unlabeled maps are its counit and the tensorators just specified. The maps $\mu_{a,b}$ preserve these half-braidings by naturality, giving the asserted coherent action on the relative center. Its effect on degrees is conjugation.

A full half-braiding yields equivariance maps
\[
 e_a:T_a(M)\xrightarrow{c_{u_a}\ot\id_{v_a}}
 (M\ot u_a)\ot v_a
 \xrightarrow{\alpha}M\ot(u_a\ot v_a)
 \xrightarrow{\id_M\ot\operatorname{coev}_a^{-1}}M.
\]
The hexagon and naturality for $m_{a,b}$ give
\begin{equation}\label{equivariant-object-law}
e_{ab}\mu_{a,b}(M)=e_aT_a(e_b),\qquad e_1=\id_M.
\end{equation}
This defines $F$ without changing the underlying object of $\D$.

Conversely, from an equivariant relative-central object $(M,c^\circ,(e_a))$, the evaluation and coevaluation isomorphisms recover 
$$c_{u_a}:u_a\ot M\to M\ot u_a.$$
For $Y\in\D_a$, use the equivalence $\D^\circ\to\D_a$, $X\mapsto X\ot u_a$, and define
\begin{equation}\label{reconstructed-halfbraiding}
\begin{split}
(X\ot u_a)\ot M&\xrightarrow{\alpha}X\ot(u_a\ot M)
 \xrightarrow{\id\ot c_{u_a}}X\ot(M\ot u_a)\\
&\xrightarrow{\alpha^{-1}}(X\ot M)\ot u_a
 \xrightarrow{c_X^\circ\ot\id}(M\ot X)\ot u_a\\
&\xrightarrow{\alpha}M\ot(X\ot u_a).
\end{split}
\end{equation}
This is natural because tensoring by $u_a$ is fully faithful. The hexagon for neutral objects is the relative hexagon. Its mixed case expresses that each $e_a$ is a morphism in the relative center, and its case $u_a,u_b$ is \eqref{equivariant-object-law}. Expanding \eqref{reconstructed-halfbraiding} for $X\ot u_a$ and $Y\ot u_b$ reduces the general homogeneous case to these identities and the ambient pentagon. Direct sums give the full half-braiding.

Restriction and reconstruction are inverse on objects and morphisms. The tensor product of half-braidings matches the equivariant tensor product defined using the tensorators of $T_a$, so $F$ is tensor. Finally, forgetting equivariance is exact, and every relative-central object $N$ is a direct summand of the underlying object of $\bigoplus_{a\in\Gamma}T_a(N)$ with its induced equivariance. Hence the functor is dominant. No averaging is used.
\end{proof}

\begin{corollary}\label{relative-exact-sequence}
There is an exact sequence of finite tensor categories, in the sense of \cite[Definition~2.7]{BN},
\[
{\rm Rep}_k(\Gamma)\longrightarrow\Z(G,\omega)
\longrightarrow\Z_{\D^{\circ}}(\D),
\] 
where the second functor is the normal dominant tensor functor of \eqref{obvious forgetful tensor functor}. 
\end{corollary}

\begin{proof}
Under Theorem \ref{gnntheorem}, the second functor forgets the action of $\Gamma$. Its kernel consists of equivariant objects, whose underlying relative-central object is a direct sum of the unit. Their equivariance is exactly a representation of $\Gamma$, giving the fully faithful first functor. For any equivariant object, the largest trivial subobject of its underlying object is preserved by the $\Gamma$-action; it is characterized as the sum of all subobjects isomorphic to the unit, and tensor autoequivalences preserve that characterization. It is therefore the underlying object of an equivariant subobject. This proves normality. Exactness and dominance were proved in Theorem \ref{gnntheorem}; together these are the asserted exact-sequence properties.
\end{proof}

For each $C\in{\rm C}(k)$, set
\begin{equation}\label{gnntheorem11}
\E_C:=\left(\bigoplus_{a\in C(k)}\Z_{\D^{\circ}}(\D_a)\right)^\Gamma.
\end{equation}
The entire sum inside the parentheses is $\Gamma$-stable. Individual components need not be stable and cannot in general be equivariantized separately.

\begin{theorem}\label{compdeco}
The equivalence of Theorem \ref{gnntheorem} restricts to equivalences of abelian categories
\[
F_C:\overline{\Z(G,\omega)_C}\xrightarrow{\simeq}\E_C.
\]
In particular, it restricts to an equivalence of tensor categories
\[
F_1:\overline{{\rm Rep}(G)_k}
\xrightarrow{\simeq}\Z(G^{\circ},\omega^{\circ})^\Gamma,
\]
where the overline denotes the Serre closure inside $\Z(G,\omega)$. The action of $\Gamma=G(k)$ on the right is the restriction of the coherent action of Theorem~\ref{gnntheorem}. In particular, its tensorators retain the mixed associator data of $\omega$.
\end{theorem}

\begin{proof}
The equivalence $F$ preserves the underlying sheaf on $G$. By the invariant-idempotent support argument, an object lies in $\overline{\Z_C}$ precisely when it is supported on $X_C$, that is, its homogeneous degrees lie in $C(k)$. Therefore, its image lies in \eqref{gnntheorem11}. Conversely, reconstruction in \eqref{reconstructed-halfbraiding} does not change the underlying sheaf of an object of $\E_C$, so that object comes from $\overline{\Z_C}$. This proves both inclusions, and hence the claimed equivalence.

For the identity class the relative center of $\D^{\circ}$ in itself is $\Z(\D^{\circ})$. Its counterpart on the left is the Serre closure of the identity-orbit category ${\rm Rep}(G)_k$, not just that smaller orbit category. Both are tensor subcategories, so the restricted equivalence is tensor.
\end{proof}

\raggedbottom

\end{document}